\documentclass[11pt,a4paper,leqno]{amsart}

\usepackage{silence}
\usepackage[T1]{fontenc}
\usepackage[utf8]{inputenc}
\usepackage[english]{babel}

\usepackage[dvipsnames]{xcolor}
\definecolor{cornellred}{rgb}{0.7, 0.11, 0.11}
\definecolor{britishracinggreen}{rgb}{0.0, 0.26, 0.15}
\definecolor{cobalt}{rgb}{0.0, 0.28, 0.67}

\usepackage[
  a4paper,
  top=2.8cm,
  bottom=2.8cm,
  left=2.8cm,
  right=2.8cm,
  marginparwidth=60pt
]{geometry}

\usepackage{amsmath}
\usepackage{amsthm}
\usepackage[utopia]{mathdesign}
\usepackage{mathtools}

\usepackage{graphicx}
\usepackage{thmtools}
\usepackage{braket}
\usepackage{caption}
\usepackage{comment}
\usepackage{enumitem}
\usepackage{multirow}
\usepackage{booktabs}
\usepackage{relsize}
\usepackage[foot]{amsaddr}
\usepackage{float}
\usepackage[colorinlistoftodos]{todonotes}
\usepackage{adjustbox}
\usepackage{rotating}
\usepackage{tikz}
\usepackage{tikz-cd}

\NewDocumentCommand{\map}{m m m O{}}{%
  \begin{tikzcd}[ampersand replacement=\&]
    #1\colon #2 \arrow[r, #4] \& #3
  \end{tikzcd}%
}
\NewDocumentCommand{\Map}{m m m O{}}{%
  \begin{tikzcd}[ampersand replacement=\&]
    #2 \arrow[r, "#1", #4] \& #3
  \end{tikzcd}%
}

\newcommand{\bigboxtimessymbol}[2]{%
  \tikz[baseline=-0.6ex, line width=#2, line cap=butt,
        line join=miter]{%
    \draw (-#1/2,-#1/2) rectangle (#1/2,#1/2);
    \draw (-#1/2,-#1/2) -- (#1/2,#1/2)
          (-#1/2, #1/2) -- (#1/2,-#1/2);
  }%
}
\DeclareRobustCommand{\bigboxtimes}{%
  \mathop{%
    \mathchoice
      {\bigboxtimessymbol{1.28em}{0.08em}}
      {\bigboxtimessymbol{0.80em}{0.08em}}
      {\bigboxtimessymbol{0.58em}{0.07em}}
      {\bigboxtimessymbol{0.43em}{0.4em}}%
  }\displaylimits
}

\usetikzlibrary{decorations.markings}
\usetikzlibrary{patterns,matrix,shapes,intersections,arrows,decorations.pathmorphing}

\usepackage[nopatch=showhyphens]{microtype}

\DeclareFontShape{T1}{mdput}{m}{up}{<->ssub*mdput/m/n}{}
\DeclareFontShape{T1}{mdput}{b}{up}{<->ssub*mdput/b/n}{}
\DeclareFontShape{T1}{mdput}{bx}{up}{<->ssub*mdput/b/n}{}

\numberwithin{equation}{section}
\allowdisplaybreaks
\usepackage[
  colorlinks,
  bookmarks,
  hypertexnames=false
]{hyperref}

\hypersetup{
  colorlinks,
  citecolor=britishracinggreen,
  filecolor=black,
  linkcolor=cobalt,
  urlcolor=purple
}

\usepackage[capitalise]{cleveref}

\DeclareSymbolFont{cmarrows}{OMS}{cmsy}{m}{n}
\SetSymbolFont{cmarrows}{bold}{OMS}{cmsy}{b}{n}
\DeclareMathSymbol{\cmminus}{\mathbin}{cmarrows}{"00}

\DeclareMathSymbol{\leftrightarrow}{\mathrel}{cmarrows}{"24}
\DeclareMathSymbol{\leftarrow}{\mathrel}{cmarrows}{"20}
   
\DeclareMathSymbol{\rightarrow}{\mathrel}{cmarrows}{"21}
   \let\to=\rightarrow
\DeclareMathSymbol{\mapstochar}{\mathrel}{cmarrows}{"37}
   \def\mapsto{\mapstochar\rightarrow}    

\DeclareSymbolFont{usualmathcal}{OMS}{cmsy}{m}{n}
\DeclareSymbolFontAlphabet{\mathcal}{usualmathcal}
\DeclareMathAlphabet\BCal{OMS}{cmsy}{b}{n} 

\DeclareMathOperator{\lHom}{\mathscr{H}\kern-0.3em\mathit{om}}
\DeclareMathOperator{\lEnd}{\mathscr{E}\kern-0.1em\mathit{nd}\,}
\DeclareMathOperator{\RRlHom}{\mathbf{R}\kern-0.025em\mathscr{H}\kern-0.3em\mathit{om}}

\newcommand{\lExt}{\mathop{\mathscr{E}\kern-0.2em\textit{xt}}\nolimits}
\DeclareMathOperator{\LTor}{{\mathscr{T}\kern-0.2em\mathit{or}}}

\newcommand{\id}{\operatorname{id}}

\newcommand{\Spec}{\operatorname{Spec}}
\newcommand{\Supp}{\operatorname{Supp}}

\newcommand{\dd}{\mathrm{d}}

\newcommand{\HH}{\mathrm{H}}
\newcommand{\OO}{\mathscr{O}}

\newcommand{\into}{\hookrightarrow}
\newcommand{\onto}{\mathrel{\mathrlap{\rightarrow}\mkern1mu\rightarrow}}
\newcommand{\derived}{\mathbf{D}}
\newcommand{\bounded}{\mathrm{b}}
\newcommand{\Db}{\derived^{\bounded}}

\makeatletter
  \newcommand{\mylabel}[2]{%
  #2%
  \phantomsection%
  \def\@currentlabel{#2}%
  \label{#1}%
  }
\makeatother

\newcommand*{\defeq}{\mathrel{\vcenter{\baselineskip0.5ex \lineskiplimit0pt
                     \hbox{\scriptsize.}\hbox{\scriptsize.}}}=}

\renewcommand{\leq}{\leqslant}
\renewcommand{\geq}{\geqslant}

\newcommand{\boldit}[1]{\boldsymbol{#1}}

\DeclareMathOperator{\MHM}{MHM}
\DeclareMathOperator{\Perv}{Perv}

\DeclareMathOperator{\Hilb}{\mathsf{Hilb}}
\DeclareMathOperator{\Quot}{\mathsf{Quot}}
\DeclareMathOperator{\Sym}{\mathsf{Sym}}

\DeclareMathOperator{\Gr}{\mathsf{Gr}}

\newcommand{\pt}{{\mathsf{pt}}}

\DeclareMathOperator{\CCoh}{\mathbf{Coh}}

\DeclareMathOperator{\Exp}{\mathsf{Exp}}
\DeclareMathOperator{\Log}{\mathsf{Log}}
\DeclareMathOperator{\Var}{Var}
\DeclareMathOperator{\Rep}{\mathsf{Rep}}
\DeclareMathOperator{\GL}{GL}

\DeclareMathOperator{\rk}{rk}

\DeclareMathOperator{\Aut}{Aut}

\newcommand{\Sn}{\Sym^{n}(\mathbb{A}^2)}

\newcommand{\Grh}{\mathsf{Gr}_w({}^{\mathfrak{p}}\mathcal{H}^k(\mathcal{F}_{n,r}))}
\newcommand{\pHcal}{{}^{\mathfrak{p}}\mathcal{H}}

\newcommand{\Fcal}{ \mathcal{F} }				

\DeclareFontFamily{OT1}{rsfs}{}
\DeclareFontShape{OT1}{rsfs}{n}{it}{<-> rsfs10}{}
\DeclareMathAlphabet{\curly}{OT1}{rsfs}{n}{it}

  

\tikzset{commutative diagrams/.cd,
mysymbol/.style={start anchor=center,end anchor=center,draw=none}}
\newcommand\MySymb[2][\square]{%
  \arrow[mysymbol]{#2}[description]{#1}}
\tikzset{
shift up/.style={
to path={([yshift=#1]\tikztostart.east) -- ([yshift=#1]\tikztotarget.west) \tikztonodes}
}
}

\theoremstyle{definition}

\newtheorem{definition}{Definition}[section]

\newtheorem{example}[definition]{Example}

\newtheorem{question}[definition]{Question}

\newtheorem{notation}[definition]{Notation}
\newtheorem{remark}[definition]{Remark}

\newtheoremstyle{thm} 
        {3mm}
        {3mm}
        {\slshape}
        {0mm}
        {\bfseries}
        {.}
        {1mm}
        {}
\theoremstyle{thm}
\newtheorem{theorem}[definition]{Theorem}
\newtheorem{corollary}[definition]{Corollary}
\newtheorem{lemma}[definition]{Lemma}
\newtheorem{proposition}[definition]{Proposition}
\newtheorem{thm}{Theorem}

\newtheorem{cor}{Corollary}[thm]

\makeatletter

\@ifundefined{c@convention}{}{}
\@ifundefined{c@thm}{}{}
\@ifundefined{c@innercustomthm}{}{}
\makeatother

\crefname{definition}{Definition}{Definitions}
\Crefname{definition}{Definition}{Definitions}
\crefname{theorem}{Theorem}{Theorems}
\Crefname{theorem}{Theorem}{Theorems}
\crefname{lemma}{Lemma}{Lemmas}
\Crefname{lemma}{Lemma}{Lemmas}
\crefname{prop}{Proposition}{Propositions}
\Crefname{prop}{Proposition}{Propositions}
\crefname{corollary}{Corollary}{Corollaries}
\Crefname{corollary}{Corollary}{Corollaries}
\crefname{remark}{Remark}{Remarks}
\Crefname{remark}{Remark}{Remarks}
\crefname{example}{Example}{Examples}
\Crefname{example}{Example}{Examples}
\crefname{assumption}{Assumption}{Assumptions}
\Crefname{assumption}{Assumption}{Assumptions}
\crefname{problem}{Problem}{Problems}
\Crefname{problem}{Problem}{Problems}
\crefname{fact}{Fact}{Facts}
\Crefname{fact}{Fact}{Facts}
\crefname{conventions}{Conventions}{Conventions}
\Crefname{conventions}{Conventions}{Conventions}
\crefname{question}{Question}{Questions}
\Crefname{question}{Question}{Questions}
\crefname{situation}{Situation}{Situations}
\Crefname{situation}{Situation}{Situations}
\crefname{openproblem}{Open Problem}{Open Problems}
\Crefname{openproblem}{Open Problem}{Open Problems}
\crefname{notation}{Notation}{Notations}
\Crefname{notation}{Notation}{Notations}
\crefname{conjecture}{Conjecture}{Conjectures}
\Crefname{conjecture}{Conjecture}{Conjectures}
\crefname{claim}{Claim}{Claims}
\Crefname{claim}{Claim}{Claims}
\crefname{observation}{Observation}{Observations}
\Crefname{observation}{Observation}{Observations}
\crefname{exercise}{Exercise}{Exercises}
\Crefname{exercise}{Exercise}{Exercises}
\crefname{convention}{Convention}{Conventions}
\Crefname{convention}{Convention}{Conventions}
\crefname{construction}{Construction}{Constructions}
\Crefname{construction}{Construction}{Constructions}
\crefname{thm}{Theorem}{Theorems}
\Crefname{thm}{Theorem}{Theorems}
\crefname{innercustomthm}{Theorem}{Theorems}
\Crefname{innercustomthm}{Theorem}{Theorems}
\crefname{cor}{Corollary}{Corollaries}
\Crefname{cor}{Corollary}{Corollaries}

\newcommand{\SE}{{\mathscr E}}

\newcommand{\SQ}{{\mathscr Q}}

\newcommand{\BA}{{\mathbb{A}}}

\newcommand{\BC}{{\mathbb{C}}}

\newcommand{\BG}{{\mathbb{G}}}

\newcommand{\BL}{{\mathbb{L}}}

\newcommand{\BN}{{\mathbb{N}}}

\newcommand{\BP}{{\mathbb{P}}}
\newcommand{\BQ}{{\mathbb{Q}}}

\newcommand{\BZ}{{\mathbb{Z}}}

\newcommand{\CA}{{\mathcal A}}

\newcommand{\CF}{{\mathcal F}}
\newcommand{\CG}{{\mathcal G}}

\newcommand{\CT}{{\mathcal T}}

\newcommand{\FS}{{\mathfrak{S}}}

\title{Purity and  Quot schemes on surfaces}
\author{Shivang Jindal, Archi Kaushik, Sergej Monavari, Andrea T. Ricolfi, Sebastian Schlegel Mejia}

\begin{document}

\begin{abstract}
Let $S$ be a smooth quasiprojective surface  and let $\SE$ be a locally free sheaf on $S$. Our main result is that the pushforward of the constant mixed Hodge module on the Quot scheme of points $\Quot^n(\SE)$ to the symmetric product $\Sym^{n}(S)$ is \emph{pure}.
Furthermore, we completely determine the decomposition theorem for the Quot-to-Chow morphism, generalizing the work of Göttsche--Soergel for Hilbert schemes. We deduce several consequences: for instance, (i) we compute the Hodge--Deligne series of $\Quot^n(\SE)$, generalizing Göttsche's formula for Hilbert schemes to higher ranks, (ii) we give a positive answer to an open question of Mozgovoy, by showing that the cohomology of $\Quot^n(\SE)$ is pure whenever the cohomology of $S$ is pure, (iii) we prove a \emph{Kirwan surjectivity} result in Borel--Moore homology for Quot schemes on surfaces, and (iv) we reprove a result of Davison, on the purity of the Borel--Moore homology of the stack of $0$-dimensional coherent sheaves on $S$.
\end{abstract}

\maketitle

{\hypersetup{linkcolor=black}
\tableofcontents}

\section{Introduction}
\subsection{Quot schemes and Hilbert schemes on surfaces} Let $S$ be a smooth quasiprojective surface over $\mathbb{C}$ and let $\SE$ be a locally free sheaf on $S$. The main goal of this article is the study and complete computation of the cohomology of the Quot scheme $\Quot^n(\SE)$,
which parametrises 0-dimensional quotients of length $n$ of $\SE$. We achieve this  by studying  the topology of the \emph{Quot-to-Chow morphism} 
\[
\Map{\pi_{n,\SE}}{\Quot^n(\SE)}{\Sym^{n}(S),}\qquad [\SE\onto \SQ]\mapsto [\Supp \SQ].
\]
When $\rk \SE=1$, the Quot scheme recovers the \emph{Hilbert scheme of points} $\Hilb^n(S)$, a smooth quasiprojective variety of dimension $2n$. In this case, the \emph{Hilbert-to-Chow morphism} $\pi_{n, \OO_S}$
is a crepant resolution of singularities of the symmetric product $\Sym^n(S)$, and  is a semismall map. 

These two key properties -- \emph{smoothness} and \emph{semismallness} -- prompted a systematic study of the cohomology of $\Hilb^n(S)$, that started  in the late 80's in the seminal works of Ellingsrud--Strømme \cite{ESHilb}, Göttsche \cite{Gott1} and Göttsche--Soergel \cite{gotsoergelhilbert} and laid the foundational basis for numerous applications to geometric representation theory, enumerative geometry and mathematical physics, see e.g.~\cite{NakHilb, Groj, OkPaQua, LehnHilb,AGT, DMVV} to mention a few.
\subsection{The cohomology of Quot schemes} In this work, we use the language of mixed Hodge modules \cite{saito} to study the topology of $\Quot^n(\SE)$ and of the Quot-to-Chow  map $\pi_{n,\SE}$. Our first main result is the following.
\begin{thm}[\Cref{main-thm-purity-A2}]\label{thm:Main thm inro}
Let $S$ be a smooth quasiprojective surface and let $\SE$ be a locally free sheaf on $S$. Denote the constant mixed Hodge module on $\Quot^n(\SE)$ by $\,\mathbb{Q}_{\Quot^n(\SE)}$. Then, the pushforward
\[ \mathcal{F}_{n,\SE}\defeq\pi_{n,\SE\ast}\mathbb{Q}_{\Quot^n(\SE)} \,\in\,\Db(\MHM(\Sym^n(S)))
\] 
is a pure mixed Hodge module of weight $0$ on the symmetric product $\Sym^n(S)$.
\end{thm}
Göttsche--Soergel's strategy \cite{gotsoergelhilbert}  for the case  $\SE=\OO_S$ is to apply the decomposition theorem to $\pi_{n, \OO_S}$, crucially relying  both on the \emph{smoothness} of $\Hilb^n(S)$ and on the  \emph{semismallness} of $\pi_{n, \OO_S}$. However, we cannot simply generalize their ideas to prove \Cref{thm:Main thm inro}. Indeed, $\Quot^n(\SE)$ is irreducible of dimension $(\rk\SE+1)n$, but \emph{singular} as soon as $\rk \SE\geq 2$ and  $n \geq 2$, and moreover the map $\pi_{n, \SE}$ is \emph{not semismall}. The singularities of $\Quot^n(\SE)$ are  notoriously hard to describe and are related to the singularities of the \emph{commuting scheme}  of $\mathfrak{gl}_{n}$. In \cite{StarkQuot} it is conjectured that $\Quot^n(\SE)$ has rational singularities.

Our next result is the computation of the generating series $ \bigoplus_{n \geq 0} \mathcal{F}_{n,\SE}$, expressed  \emph{plethystically}. For a commutative $\BN$-graded monoid $(X,\cup)$ in the category of $\BC$-schemes, consider the category $\widehat \derived(X) = \prod_{n\geq 0}\Db(\MHM(X_n))$. See \Cref{subsec: Plethystic operations} for the definition of the functor $\Sym_{\cup}$, landing in $\widehat \derived(X)$.
\begin{thm}
[\cref{thm:sheaf and coho as Sym} and \cref{prop:computation of Hn}]
\label{thm: Computation of decomposition}
Let $S$ be a smooth quasiprojective surface and let $\SE$ be a rank $r$ locally free sheaf on $S$. 
Let $q_{(m)}\colon S\into \Sym^{m}(S)$ be the inclusion of the small diagonal. Set $\mathbb{L}\defeq  \mathrm{H}_c(\mathbb{A}^1,\mathbb{Q})$. Then, there is an isomorphism
\begin{equation}
\label{eq:FV-is-Sym}
\bigoplus_{n \geq 0} \mathcal{F}_{n,\SE} \simeq \Sym_{\cup}\left(\bigoplus_{m \geq 1} q_{(m)\ast}\mathbb{Q}_{S}\otimes \HH^{\ast} (\mathbb{P}^{r-1}) \otimes \mathbb{L}^{(m-1)r}\right)\,\in\,\widehat \derived(\Sym(S))
\end{equation}
of mixed Hodge modules on $\Sym(S)$. In particular, $\CF_{n, \SE}$ only depends on $\rk\SE$. 
\end{thm}
In the case $\SE=\OO_S$, \Cref{thm: Computation of decomposition} recovers the decomposition of \cite{gotsoergelhilbert}, see \Cref{example: Thm B generalizes Gottsche Soeergel}. 

Given the purity obtained in \Cref{thm:Main thm inro}, and using the fact that the category of pure mixed Hodge modules is semisimple by Saito's Decomposition theorem \cite{saito}, we deduce that $\mathcal{F}_{n,\SE}$ admits a decomposition into shifts of simple mixed Hodge modules whose supports are strata of a natural stratification of $\Sym^n(S)$, described in \Cref{sec:sym-strata}. The proof of \Cref{thm: Computation of decomposition} follows by the explicit computation of the simple pieces in this decomposition of $\CF_{n, \SE}$, suitably rearranged with the plethystic formalism.

Pushing forward \eqref{eq:FV-is-Sym} to a point, we compute  the mixed Hodge structure of the (compactly supported) cohomology  of $\Quot^n(\SE)$. In what follows, we denote by $\HH^{\ast}_{(c)}(-)$ either the \emph{cohomology } or the \emph{compactly supported cohomology} functor.

\begin{cor}[\Cref{thm:sheaf and coho as Sym}]
\label{cor: coho int for quot}
Let $S$ be a smooth quasiprojective surface and let $\SE$ be a rank $r$ locally free sheaf on $S$. 
Then, there is an isomorphism of $\BZ$-graded mixed Hodge structures
\begin{equation}
\label{eqn:Coh-as-Sym}
\bigoplus_{n \geq 0} \HH^{\ast}_{(c)}(\Quot^n(\SE)) \simeq \Sym \left( \bigoplus_{m \geq 1} \HH^{\ast}_{(c)}(S) \otimes \HH^{\ast} (\mathbb{P}^{r-1}) \otimes \mathbb{L}^{(m-1)r}\right)\,\in\,\widehat \derived(\BN).
\end{equation}
\end{cor}
Our next result uses \Cref{cor: coho int for quot} to compute cohomological invariants of $\Quot^n(\SE)$. Let
\begin{align*}
\mathsf h_{(c)}(S;u,v,z) = \sum_{p,q,k}h_{p,q,{(c)}}^ku^pv^qz^k
\end{align*}
be the \emph{Hodge--Deligne series} of $S$,  where  $h_{p,q,{(c)}}^k = h^{p,q}(\HH^k_{(c)}(S))$, and let
\begin{align*}
        \mathsf{HD}^{(c)}_r(t)=\sum_{n\geq 0}\mathsf h_{(c)}(\Quot^n(\SE);u,v,z)t^n
\end{align*}
be the generating series of the Hodge--Deligne series of $\Quot^n(\SE)$, see \eqref{eqn: HD graaaa}. We deduce an infinite product formula for $\mathsf{HD}^{(c)}_r(t)$ as well as for the Poincaré series of $\Quot^n(\SE)$. See \Cref{sec:HD-series} for the definition of the operator $\Sym$ acting on $\BZ[u^{\pm 1},v^{\pm 1},z^{\pm 1}]$.

The following result generalizes Göttsche's formula \cite{Gott2} from   Hilbert schemes to Quot schemes.

\begin{thm}[\Cref{thm: HD poly}, \Cref{cor: Poincare}]\label{thm: intro HD}
Let $S$ be a smooth quasiprojective surface and let $\SE$ be a locally free sheaf of rank $r$ on $S$. Set $\mathsf q=\mathsf h_c(\BA^1;x,y,z) = xyz^2$. There are identities
\begin{align*}
  \mathsf{HD}^{(c)}_r(t)&= \Sym \left(\mathsf h_{(c)}(S;u,v,z)\frac{1-\mathsf q^r}{1-\mathsf q}\frac{t}{1-t\mathsf q^r}\right) \\
&=\prod_{m\geq 1}\prod_{j=0}^{ r-1}\prod_{p,q,k}\left(1-(-1)^kt^mu^{p+j+r(m-1)} v^{q+j+r(m-1)}z^{k+2j+2r(m-1)}\right)^{(-1)^{k+1}h_{p,q,(c)}^k}.    
\end{align*}
At the level of Poincar\'e polynomials, we have the identities
\begin{equation*}
\label{eq:Poincare poly intro}
\sum_{n\geq 0}P_{(c)}(\Quot^n(\SE),z)t^n=
\prod_{m\geq 1}\prod_{j=0}^{r-1}\prod_{i=0}^4 \left(1-(-1)^it^m z^{2j+2r(m-1)+i}\right)^{(-1)^{i+1}b_i^{(c)}}.
\end{equation*}
\end{thm}
 We remark that the only known computation of  the Poincaré polynomials of $\Quot^n(\SE)$ (prior to our work), when the Quot scheme is singular, was done for $S=\BA^2$ and $n=2$ in \cite{Pielasa}.
 
\smallbreak
 As a corollary of the explicit formulas computed in \Cref{thm: intro HD}, we prove a \emph{rank $r$-to-$1$} factorization of $\mathsf{HD}^{(c)}_r(t)$, phrased in terms  of the natural  \emph{power structure} on $\BZ[u^{\pm 1},v^{\pm 1},z^{\pm 1}]$.
\begin{cor}[\Cref{prop: rank r to 1}]
\label{cor:factorization-intro}
Let $S$ be a smooth quasiprojective surface and let $\SE$ be a locally free sheaf of  rank $r$ on $S$. Set $\mathsf q= xyz^2$. Then 
\[
\mathsf{HD}^{(c)}_r(t) = \prod_{j=0}^{r-1} \mathsf{HD}^{(c)}_1(t\mathsf q^{r-1})^{\mathsf q^{j+1-r}}.
\]
\end{cor}
A similar rank $r$-to-$1$ factorization phenomenon appears for motives of Quot schemes on smooth quasiprojective curves \cite{Monavari-Ricolfi-motive-nested-quot, monavari2024hyperquot,Ricolfi-Virtual-Classes-Quot} and for \emph{virtual invariants} of more general Quot schemes \cite{Fasola-Monavari-Ricolfi,FM_tetra, monavari2026gaugeorigamibrokenlines}. In all these cases, the factorization follows as an application of (virtual) torus localization. 

Let us also note that \Cref{cor:factorization-intro} suggests a \emph{categorification} statement, namely the existence of a semiorthogonal decomposition of $\Db(\Quot^n(\SE))$ in terms of components of the form $\Db\left(\prod_i \Hilb^{n_i}(S)\right)$ indexed by the condition $n_1+\cdots+n_r=n$, possibly after passing to derived enhancements of Quot and Hilbert schemes. We leave this for future research.

\subsection{Mozgovoy's question on purity} In \cite{mozgovoy2019motivic}, Mozgovoy derives a formula for the motive (i.e.~the class in the Grothendieck ring of varieties) of $\Quot^n(\SE)$. The main result in \cite{mozgovoy2019motivic} states that
\begin{equation}
\label{eqn:Mozgovoy}
\sum_{n\geq 0}\,[\Quot^n(\SE)]t^n = \Exp\left(\frac{[S\times \BP^{r-1}]t}{1-\mathbf L^rt} \right)\,\in\,K_0(\Var_{\BC})\llbracket t \rrbracket.
\end{equation}
 For $S=\mathbb{A}^{2}$ and $S=\mathbb{P}^{2}$, the coefficients in the expression \eqref{eqn:Mozgovoy} are polynomials in the Lefschetz class $\mathbf{L} = [\BA^1]$. This suggests the existence of affine pavings that would yield purity for the compactly supported cohomology of the Quot schemes on $\BA^2$ and $\BP^2$. Therefore, Mozgovoy asks the following question.

\begin{question}[Mozgovoy {\cite{mozgovoy2019motivic}}]\label{Question: Mozgovoy}
Let $S$ be either $\BA^2$ or a smooth projective surface and let be $\SE$ a locally free sheaf on $S$. Is the mixed Hodge structure on the cohomology $\HH_c^*(\Quot^n(\SE))$ pure?
\end{question}

Since the $\Sym$ operation preserves purity, \Cref{cor: coho int for quot} answers this question in the positive, as the following (more general) statement illustrates.

\begin{thm}[\Cref{cor: Purity of Quot}]
\label{cor:Ans to Mozgovoy question intro}
  Let $S$ be a smooth quasiprojective surface such that $\HH^{\ast}_{(c)}(S)$ is a pure mixed Hodge structure of weight $0$. Then, for all $n\geq 0$, we have that $\HH^{\ast}_{(c)}(\Quot^n(\SE))$ is a pure mixed Hodge structure of weight $0$, i.e.~$\HH^k_{(c)}(\Quot^n(\SE))$ is a pure mixed Hodge structure of weight $k$ for all $k \geq 0$.  
\end{thm}

Given the computation of the Poincaré series  of $\HH^{*}(\Quot^n(\SE))$ realized in \Cref{thm: intro HD}, following Nakajima and Grojnowski's \cite{NakHilb, Groj} categorification of Göttsche's formula  and the recent works \cite{marian2026derivedcategoriesquotschemes,MarianNegutcoho,jindallim,kaushik2026cohomologyhyperquotschemescurves} for Quot schemes of curves, 
it is natural to ask whether a PBW basis of the positive-half of  a shifted Yangian would yield a basis for $\HH^{*}(\Quot^n(\SE))$. 
\subsection{Kirwan surjectivity}\label{subsection: KIrwan Surjectivity intro}
Denote by  $\CCoh^{n}(S)$ the stack of 0-dimensional coherent sheaves on $S$ of length $n$, and consider the stack $\mathrm{Fr}_{n}(\SE)$ parametrizing pairs $(\SQ,\phi)$, where $\SQ$ is a 0-dimensional coherent sheaf of length $n$ on $S$ and $\phi\colon \SE\rightarrow \SQ$ is a morphism. The natural forgetful map
\begin{equation}
\label{eq:map from Fr to Coh}
\begin{tikzcd}
\mathrm{Fr}_{n}(\SE)\arrow{r}{p} & \CCoh^{n}(S),\qquad (\SQ,\phi) \mapsto \SQ 
\end{tikzcd}
\end{equation}
realizes $\mathrm{Fr}_{n}(\SE)$ as a rank $rn$ vector bundle over the stack $\CCoh^{n}(S)$. 
The Quot scheme $\Quot^n(\SE)$ is the open substack of $\mathrm{Fr}_{n}(\SE)$ consisting of those pairs $(\SQ,\phi)$ such that $\phi$ is surjective. In fact, the space $\Quot^n(\SE)$ is the semistable locus of a natural $\Theta$-stratification \cite{zbMATH07272594} on the stack $\mathrm{Fr}_{n}(\SE)$, see \cite{Todaquot} for the similar case of  Quot schemes on curves.

In similar settings, the \emph{Kirwan surjectivity} in Borel--Moore homology holds, see e.g.~\cite{arXiv:2509.21298,zbMATH07272594}, provided that the (derived) analogue of the stack $\mathrm{Fr}_{n}(\SE)$ satisfies a certain technical condition. However, as explained in \Cref{rem:Kirwan surjectivity lit}, in our case this technical condition does not hold. Nevertheless, we can prove Kirwan surjectivity for the Quot scheme $\Quot^n(\SE)$ as a consequence of \cref{thm:Main thm inro}.

\begin{thm}[\Cref{cor:Kirwan Surjectivity}]\label{Cor: Kirwan Surjectivity intro}
Let $S$ be a smooth quasiprojective surface and let $\SE$ be a locally free sheaf of rank $r$ on $S$. Then the open pullback in Borel--Moore homology induces a surjection of mixed Hodge structures
\[
\begin{tikzcd}
\HH_{\ast}^{\mathrm{BM}}(\CCoh^{n}(S))\otimes\mathbb{L}^{-rn}\cong \HH_{\ast}^{\mathrm{BM}}(\mathrm{Fr}_{n}(\SE))\arrow[two heads]{r} & \HH_{\ast}^{\mathrm{BM}}(\Quot^n(\SE)).
\end{tikzcd}
\]
\end{thm}

\subsubsection{Borel--Moore homology of $\CCoh^{n}(S)$}
Our final application of \Cref{thm:Main thm inro} is an independent proof of the purity of the Borel--Moore homology of $ \CCoh^{n}(S)$ (already proved in \cite{Davisonpreprojective}, see also \cite{BDHSM}) via its interpretation  as the \emph{stabilization} of the pure mixed Hodge structure on the Borel--Moore homology of $\Quot^n(\OO^{\oplus r}_S)$.
\begin{cor}[\Cref{prop: Homology of Quot stabilizes to that of Coh}, \Cref{cor:Purity of Coh}]
Let $S$ be a smooth quasiprojective surface. Then, there is an isomorphism
\[
\HH_{*}^{\mathrm{BM}}(\CCoh^{n}(S))\cong \varprojlim_r\HH_{*}^{\mathrm{BM}}(\Quot^n(\OO^{\oplus r}_S))\otimes \mathbb{L}^{rn}.
\]
In particular, if $\,\HH_{*}^{\mathrm{BM}}(S)$ is pure of weight $0$, the Borel--Moore homology $\HH_{*}^{\mathrm{BM}}(\CCoh^{n}(S))$ is pure of weight $0$.
\end{cor}

\subsection{Strategy of the proof of \Cref{thm:Main thm inro}}
To prove our first main result, we first observe that by the \'etale-local nature of the problem (see \Cref{lemma:purity-checked-on-etale-cover}), it is enough to treat the case $S=\BA^2$. Purity of the mixed Hodge module $\CF_{n,r}\defeq \pi_{n,\OO^{\oplus r}\ast}\mathbb{Q}_{\Quot^n(\OO_{\BA^2}^{\oplus r})}$ is equivalent to the vanishing
\begin{equation}
\label{eq:required vanishing intro}
\Grh=0, \qquad \mbox{for all } w\neq k\in\mathbb{Z}.
\end{equation}
The $\BG_m$-action $t\cdot (x,y)=(tx, t^{-1}y)$ lifts to $\Sym^n(\BA^2)$. The attracting locus is $g_n\colon \Sym^n(\BA^1 \times \set{0})\into \Sym^n(\BA^2)$, and the fixed locus is $f_n\colon n[0] \into \Sym^n(\BA^1 \times \set{0})$. We define the \emph{hyperbolic localization} functor (see~\Cref{sec: hyperbloc}) by
\[
\Phi_n(\CF)=f_n^{!}g_n^*\mathcal{F}\in \Db(\mathrm{MHS}),
\]
analogously to Braden--Mautner \cite{bradenmautner}. We show the key vanishing 
\begin{equation*}
\Phi_n\left(\Grh\right)=0, \qquad \mbox{for all } w\neq k\in\mathbb{Z},
\end{equation*}
exploiting the properties of $\Phi_n$ and the purity of the compactly supported cohomology of \emph{semi-nilpotent quiver varieties} proved in \cite{schiffmannvasserot} (cf.~\Cref{purityondiagonal}).
Then, \Cref{onthediagonal} reduces the problem to checking that $\Grh$ is supported on the small diagonal $\BA^2 \into \Sym^n(\BA^2)$, which we confirm thanks to a factorization property of $\CF_{n,r}$, cf.~\Cref{prop:sheaffactorization}.

\subsection*{Conventions}
We work over the field $\BC$ of complex numbers. A \emph{variety} is a quasiprojective reduced scheme over $\BC$. If $Y$ is a variety, we set $\Sym^n(Y) = Y^n / \mathfrak S_n$ and $\Sym (Y) = \coprod_{n\geq 0}\Sym^n(Y)$. 
Sheaves and cohomology groups have rational coefficients. Let $Y$ be a variety and let $H$ be a rational vector space or a mixed Hodge structure. We denote the associated constant local system or variation of mixed Hodge structures on $Y$ using the same notation $H$, unless stated otherwise. For compactness of notation, we denote by $\HH^{\ast}_{(c)}(-)$ either the cohomology or the compactly supported cohomology functor. Partitions $\lambda \vdash n$ of an integer $n$ are written as weakly decreasing sequences $(\lambda_1,\ldots,\lambda_{l(\lambda)})$, with $\sum_i\lambda_i = n$, and we set $G_\lambda = \Aut(\lambda)$.

\subsection*{AI disclosure}
During the preparation of this article, the authors used  ChatGPT 6.2 Astra and Claude Opus 5.5 to assist with literature search and mathematical verifications. The mathematical ideas and proofs are the work of the authors, who have written and independently verified the entire manuscript. The authors assume full responsibility for the mathematical content and conclusions of this work.

\subsection*{Acknowledgements}
We are grateful to Mathis Duguin, Camilla Felisetti,  Niccolò Giacomini, Nicolas Hemelsoet, Lucien Hennecart, Andres Fernandez Herrero, Woonam Lim, Elsa Maneval, Andrei Neguț and Andrés Ibáñez Núñez, for useful discussions. S.M. is supported  by  the HORIZON-MSCA-2024-PF-01 Project 
101203281 ``Moduli Spaces of Sheaves: Geometry and Invariants'', funded by the Research and Innovation framework programme Horizon Europe. S.S.M.\ was supported by the Swiss National Science Foundation [no. 218340]. S.M. and A.R. are members of GNSAGA-INDAM.

\section{Mixed Hodge modules}\label{section:MHMs}
In this section, following \cite{saito,saitoguide,benpurity, Tubach2025}, we set up the notations regarding mixed Hodge modules and recall some of their important properties needed throughout the paper. We work with $\mathbb{Q}$-coefficients. All functors are appropriately derived.

\subsection{Review of mixed Hodge modules and purity}\label{sec: MHM}
Let $X$ be a separated finite type scheme over $\mathbb{C}$. We denote by $\MHM(X)$ Saito's abelian category of \emph{mixed Hodge modules} on $X$ and by $\Db(\MHM(X))$ its bounded derived category. We denote by $\Perv(X)$ the abelian category of \emph{perverse sheaves} on $X$, namely the heart of the t-structure on the derived category $\Db_c(X)$ of constructible sheaves determined by the middle-perversity function. 
There is an exact faithful functor 
\[
\map{\mathrm{rat}_X}{\MHM(X)}{\Perv(X),}
\]
which extends to a conservative functor of derived categories
\[
\map{\mathrm{rat}_X}{\Db(\MHM(X))}{\Db_c(X),}
\]
compatible with the six-functor formalism. In particular, the natural heart of the category $\Db(\MHM(X))$ corresponds to  the perverse heart of $\Db_c(X)$. We denote the cohomology functor with respect to these natural t-structures by $\pHcal$. We have 
\[ 
\mathrm{rat}_X(\pHcal(\mathcal{F})) = \pHcal(\mathrm{rat}_X(\mathcal{F})).
\] 
Over a point, $\MHM(\pt)$ is the category $\mathrm{MHS}$ of graded-polarizable mixed Hodge structures. In general, a mixed Hodge module $\mathcal{F} \in \MHM(X)$ comes with a finite increasing weight filtration 
 \[ 
 W_{\bullet}(\mathcal{F})\colon \quad \cdots \subset W_w(\mathcal{F}) \subset W_{w+1}(\mathcal{F}) \subset \cdots
 \]
with graded pieces $\Gr_{w}(\mathcal{F})\defeq  W_{w}(\mathcal{F})/W_{w-1}(\mathcal{F})$. A mixed Hodge module $\mathcal{F}$ is \emph{pure of weight} $w$ if $\Gr_k(\mathcal{F})=0$ for all $k \neq w$. Any morphism of mixed Hodge modules strictly preserves the weight filtration and thus a morphism between pure mixed Hodge modules of distinct weights is zero.

A complex $\mathcal{F} \in \Db(\MHM(X))$ is said to have weights $\leq a$ if, for all $k \in \mathbb{Z}$, one has
\[ \Gr_w(\pHcal^{k}(\Fcal))=0\quad \text{ for all}\ w>a+k.
\] 
and similarly it is said to have weights $\geq a$ if, for all $k \in \mathbb{Z}$, one has
\[ \Gr_w(\pHcal^{k}(\Fcal))=0\quad \text{ for all}  \ w<a+k.
\] 
It is said to be \emph{pure of weight} $a$ if both conditions hold. In particular, a complex of mixed Hodge modules $\mathcal{F}$ is pure of weight $0$ if and only if $\pHcal^k(\mathcal{F})$ is pure of weight $k$ for all $k$. We denote by 
\[
\mathbb{L}\defeq  \mathrm{H}_c(\mathbb{A}^1,\mathbb{Q}) 
\] 
the complex of pure Hodge structures, concentrated in cohomological degree $2$. 

For any morphism $f\colon  X\rightarrow Y$, the derived functors $f^*,f_{!}$ can only decrease weights and $f_{*}, f^!$ can only increase weights. Thus, since $f_\ast = f_!$ for a proper morphism $f$, pushforward along proper morphisms preserves purity, cf.~\cite[p.~324]{saito}.

If $g\colon U \rightarrow X$ is an \'etale morphism, then since $g^{!}=g^*$, we have that $g^*$ preserves weights. Additionally, $g^*$ commutes with the perverse t-structure. Thus, for $\mathcal{F} \in \Db(\MHM(X))$, we have \begin{equation*}
g^*\Gr_w(\pHcal(\mathcal{F})) \simeq \Gr_{w} (\pHcal(g^*\mathcal{F})).  
\end{equation*} 
Furthermore, we can use \'etale covers to detect purity of mixed Hodge modules.
\begin{lemma}[{\cite[Lemma 2.2]{benpurity}}]
\label{lemma:purity-checked-on-etale-cover}
Let $\{g_i\colon U_i \rightarrow X\}_{i \in I}$ be an \'etale cover of $X$. Then a complex $\mathcal{F} \in \Db(\MHM(X))$ is pure of weight $0$ if and only if the complex $g_i^*\mathcal{F}$ is pure of weight $0$ for all $i \in I$. 
\end{lemma}
\subsubsection{Intermediate extension functor and intersection complex}
We recall some properties of the intermediate extension functor from \cite{saitonotes}. Let $j\colon U \into X$ be a locally closed immersion. Then we define the intermediate extension functor 
\[ 
\begin{tikzcd}
\Db(\MHM(U)) \arrow{r}{j_{!\ast}} & \MHM(X)
\end{tikzcd}
\] 
as $j_{!\ast}(\mathcal{F}) \defeq  \mathrm{Im}(\pHcal^0(j_! \mathcal{F}) \rightarrow \pHcal^0(j_{\ast}\mathcal{F}))$. Let $i\colon Z \into X$ be any integral closed subscheme and let $k\colon U \into Z$ be a smooth dense open subset. Let $\mathbb{V}$ be a polarizable variation of mixed Hodge structures on $U$. Then we define the \emph{intersection complex} as 
\[
\mathrm{IC}_Z(\mathbb{V})\defeq  (ik)_{!\ast}(\mathbb{V}[\dim Z]).
\] 
Note that, since $i$ is a closed immersion, by the definition of the intermediate extension functor, we have that $\mathrm{IC}_Z(\mathbb{V}) \simeq i_* k_{!*}(\mathbb{V}[\dim Z])$. If $X$ is integral, then $\mathrm{IC}_X(\mathbb{V})$ is the unique mixed Hodge module on $X$, supported on $Z$, which restricts to $\mathbb{V}[\dim Z]$ on the dense open $U$ and has no nonzero subobject or nonzero quotient supported on $Z\backslash U$.

\subsection{Stratifications}
Following \cite[Sec.~2.3]{achar}, 
let $\mathcal{T}$ be a stratification of $X$ and let $\Db_{\mathcal{T}}(X)$ be the full
subcategory of $\Db_c(X)$ consisting of complexes $\mathcal{F}$ whose
cohomology sheaves $\mathcal{H}^k(\mathcal{F})$, with respect to the standard heart, restrict to local systems on
each stratum, i.e.~for each stratum $T\in \CT$, the constructible sheaf $\mathcal{H}^k(\Fcal)|_{T}$ is a local system. We can similarly define the category of $\mathcal{T}$-constructible perverse sheaves
$\Perv_{\mathcal{T}}(X)\subset \Perv(X)$. Let us define the subcategories
\begin{align*}
\MHM_{\mathcal{T}}(X)
&=\set{\Fcal \in\MHM(X)\mid
\mathrm{rat}_X(\Fcal)\in\Perv_{\mathcal{T}}(X)}\subset \MHM(X),\\
\Db_{\mathcal{T}}(\MHM(X))
&=\set{\mathcal{F}\in\Db(\MHM(X))\mid
\mathrm{rat}_X(\Fcal)\in\Db_{\mathcal{T}}(X)}\subset \Db(\MHM(X)).
\end{align*}

By definition of the category of perverse sheaves, the perverse truncation functor and hence the perverse cohomology functors preserve $\mathcal{T}$-constructibility. Moreover, $\mathcal{T}$-construcible perverse sheaves are closed under taking subobjects and quotients. Thus, for an object $\mathcal{F} \in \Db_{\mathcal{T}}(\MHM(X))$, we have that $\Gr_w(\pHcal^k(\mathcal{F})) \in \MHM_{\mathcal{T}}(X)$.

\subsection{Plethystic operations}
\label{subsec: Plethystic operations}
Let $(X,\cup)$, with $X=\coprod_{d\geq 0}X_d$, be an $\BN$-graded commutative monoid in the category of $\BC$-schemes, such that $X_0 = \Spec \BC$. If $\lambda = (\lambda_1,\ldots,\lambda_l)$ is a sequence of positive integers (e.g.~a partition) summing up to $d$, let $\cup_{\lambda} \colon X_{\lambda_1}\times \cdots \times X_{\lambda_l} \to X_d$ be the associated monoid operation. Define 
\[
\widehat{\derived}(X) = \prod_{d\geq 0}\Db(\MHM(X_d)), \qquad \widehat{\derived}_{>0}(X) = \prod_{d>0}\Db(\MHM(X_d)).
\]
We have a product
\[
\Map{\boxtimes_\cup}{\widehat{\derived}(X)\times \widehat{\derived}(X)}{\widehat{\derived}(X),}\qquad 
(\CF,\CG) \mapsto \CF \boxtimes_\cup \CG = \cup_{!} (\CF \boxtimes \CG).
\]
Explicitly, if $\CF$ and $\CG$ have components $\CF_d$ and $\CG_d$ respectively, the $d$-th component of $\CF \boxtimes_\cup \CG$ is given by
\[
(\CF \boxtimes_\cup \CG)_d = \bigoplus_{a+b=d}\cup_{a,b!}(\CF_a \boxtimes \CG_b) \in\Db(\MHM(X_d)).
\]
Assume $\cup \colon X \times X \to X$ is a finite morphism. By (the proof of) \cite[Prop.~3.5]{davisonmeinhardt}, one has that $\boxtimes_\cup$ is symmetric monoidal, exact and preserves weight filtrations.

Given an object $\CF \in \widehat{\derived}_{>0}(X)$, with components $(\CF_d)_{d>0}$, we produce a symmetrized object $\Sym_\cup(\CF) \in \widehat{\derived}(X)$ as follows. Let $\cup^l \colon X^l \to X$ be the $l$-th iteration of $\cup$. We set 
\[
\Sym^l_\cup(\CF) \defeq (\CF^{\boxtimes_\cup l})^{\FS_l} = (\cup^l_! \CF^{\boxtimes l})^{\FS_l}\,\in\,\widehat{\derived}(X).
\]
The $d$-th component of $\Sym^l_\cup(\CF)$ is given by
\[
\Sym^l_\cup(\CF)_d = \bigoplus_{\substack{\lambda \vdash d \\ l(\lambda)=l}} \left(\CF_{\lambda_1}\boxtimes_\cup\cdots\boxtimes_\cup \CF_{\lambda_l}\right)^{G_\lambda} = \bigoplus_{\substack{\lambda \vdash d \\ l(\lambda)=l}} \left(\cup_{\lambda!}\bigboxtimes_{i=1}^l \CF_{\lambda_i}\right)^{G_\lambda}
\]
where the sum is over partitions $\lambda = (\lambda_1,\ldots,\lambda_l)\vdash d$ and $G_\lambda = \Aut(\lambda)$. Note that $\Sym^l_\cup(\CF)_d = 0$ when $l>d$. This construction defines a functor
\[
\map{\Sym_\cup}{\widehat{\derived}_{>0}(X)}{\widehat{\derived}(X),} \qquad \CF \mapsto \Sym_\cup(\CF) \defeq \bigoplus_{l\geq 0}\Sym^l_\cup(\CF),
\]
preserving purity, cf.~\cite[Prop.~3.8]{davisonmeinhardt}. Here the direct sum is taken component-wise, i.e.~the $d$-th component of $\Sym_\cup(\CF)$ is
\[
\Sym_\cup(\CF)_d = \bigoplus_{l=0}^d\Sym_\cup^l(\CF)_d = \bigoplus_{\lambda\vdash d}\left(\CF_{\lambda_1}\boxtimes_\cup\cdots\boxtimes_\cup \CF_{\lambda_r}\right)^{G_\lambda}\in\Db(\MHM(X_d)).
\]
 We will mainly focus on the following two special cases.
\begin{enumerate}
    \item [\mylabel{Sym-N}{(1)}] Set $(X,\cup)=(\BN,+)$. In this case, objects of $\widehat{\derived}(\BN)$ are sequences of objects in $\Db(\mathrm{MHS}) = \Db(\MHM(\pt))$, and we simply write $\Sym \defeq \Sym_+$.
    \item [\mylabel{Sym-Sym}{(2)}] Set $(X,\cup)=(\Sym(Y),\cup)$ for a variety $Y$. Here $\cup \colon \Sym (Y) \times \Sym (Y) \to \Sym (Y)$ is the finite map given by addition of cycles.
\end{enumerate}

\section{Quot schemes and symmetric products}\label{Section:Quot and sym}
\subsection{The Quot scheme of points}
Let $S$ be a smooth quasiprojective complex surface, and let $\SE$ be a locally free sheaf of rank $r$. Let $n \in \mathbb Z_{\geq 0}$ be an integer. The \emph{Quot scheme of $n$ points} attached to $(S,\SE,n)$ is the fine moduli space $\Quot^n(\SE)$ parametrising isomorphism classes of quotients $\SE \onto \SQ$, where $\SQ$ is a 0-dimensional sheaf of length $\chi(\SQ) = n$. It is known to be irreducible of dimension $n(r+1)$, and singular for $n \geq 2$ as long as $r\geq 2$, cf.~\cite{Ellingsrud_Lehn} and \cite[Ex.~3.3]{Cazzaniga-Ricolfi-Framedsheaves}. There is a smooth morphism
\begin{equation}
\label{eqn:quot-to-coh}
\map{\theta_{n,\SE}}{\Quot^n(\SE)}{\CCoh^{n}(S)}
\end{equation}
of relative dimension $rn$ to the stack $\CCoh^{n}(S)$ of 0-dimensional sheaves of length $n$ over $S$, cf.~\cite[Lemma 1.24]{Fantechi-Ricolfi-structural} or \cite[Cor.~1.6]{Ricolfi-Savvas}. Furthermore, there is a `support morphism'
\[
\begin{tikzcd}
\sigma_n\colon \CCoh^{n}(S)\arrow{r} & \Sym^n(S)
\end{tikzcd}
\]
to the symmetric product, see \cite{Rydh1,Fantechi-Ricolfi-structural,Fantechi-Ricolfi-motivic}. It coincides with the good moduli space morphism for the algebraic stack $\CCoh^{n}(S)$. The composition
\[
\map{\pi_{n,\SE}\defeq \sigma_n\circ \theta_{n,\SE}}{\Quot^n(\SE)}{\Sym^n(S),}
\]
is a proper morphism called the \emph{Quot-to-Chow morphism}; it takes a closed point $[\SE \onto \SQ]$ to the effective 0-cycle $[\Supp \SQ]$.

\begin{notation}
\label{def:punctual-quot}
In the local case $S=\BA^2$, we set $\Quot_{n,r} \defeq \Quot^n(\OO_{\BA^2}^{\oplus r})$, and we denote by $\pi_{n,r}$ the Quot-to-Chow map $\Quot_{n,r} \to \Sym^n(\BA^2)$. The \emph{punctual Quot scheme} is the fibre of $\pi_{n,r}$ over the fat cycle $n[0]\in \Sym^n(\BA^2)$. It is denoted $\Quot_{n,r}^0$. 
\end{notation}

By a standard argument, cf.~for instance \cite[Sec.~2.1]{Ricolfi-Power-structure-Quot}, $\Quot_{n,r}^0$ is (non-canonically) isomorphic to $\pi_{n,\SE}^{-1}(n[p])$ for any smooth surface $S$, any  rank $r$ bundle $\SE$ and any point $p \in S$.

\subsection{Quiver models for the local Quot scheme}
\label{sec:quot-via-quivers}
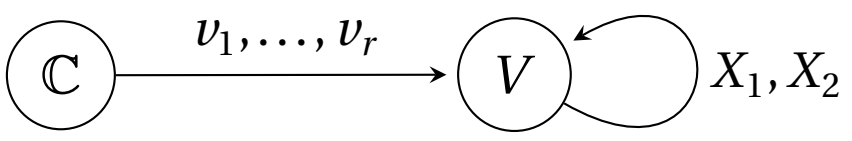
\begin{figure}[H]
\centering
\begin{tikzpicture}[>=stealth,shorten >=2pt,looseness=.5,auto]
  \matrix [matrix of math nodes,
           column sep={3cm,between origins},
           row sep={3cm,between origins},
           nodes={circle, draw, minimum size=5.5mm}]
{ 
|(A)| \BC & |(B)| V \\         
};
\tikzstyle{every node}=[font=\small\itshape] 
\node [anchor=west,right] at (2.65,0) {$X_1,X_2$}; 
\draw[->] (A) -- node[above] {$v_1,\ldots,v_r$} (B);
\draw[->] (B) to [out=-30,in=30,looseness=9] (B) node {};
\end{tikzpicture}
\caption{A $(1,n)$-dimensional representation of the $r$-framed 2-loop quiver $\boldit{Q}_{2}^r$.}
\label{fig:ADHM-quiver}
\end{figure}

In the local case $S=\BA^2$, one realizes $\Quot_{n,r}$ as a GIT quotient as follows. Consider an $n$-dimensional complex vector space $V$, and let $R\cong \BA^{2n^2+rn}$ be the representation space of the quiver $\boldit{Q}_{2}^r$ depicted in \Cref{fig:ADHM-quiver}. 
Consider the open subscheme $U \subset R$ $(X_1,X_2,v_1,\ldots,v_r) \in R$ of \emph{cyclic} tuples: this means that the vectors $v_1,\ldots,v_r \in V$ generate the unframed representation $(X_1,X_2) \in \Rep_n(\boldit{L}_{2})$, where $\boldit{L}_{2} \subset \boldit{Q}_{2}^r$ is the 2-loop quiver; equivalently, the smallest subspace of $V$ preserved by $X_1,X_2$ and containing $v_1,\ldots,v_r$ is $V$ itself. Then $U$ coincides with the stable locus with respect to a King stability condition on $\boldit{Q}_{2}^r$ (see e.g.~the proof of \cite[Lemma 2.2]{Cazzaniga-Ralaivaosaona-Ricolfi}), and $\GL_n$ acts freely on $U$. The quotient $U/\GL_n$ is a smooth quasiprojective variety of dimension $n^2+rn$, called the \emph{noncommutative Quot scheme}. Let $Z \subset U$ be the closed subscheme cut out by the commutator relation $[X_1,X_2] = 0$. Then we have a commutative diagram
\[
\begin{tikzcd}
Z\MySymb{dr}\arrow[hook]{r}\arrow{d} & U\arrow[hook]{r}\arrow{d} & R \\
\Quot_{n,r}\arrow[hook]{r} & U/\GL_n & 
\end{tikzcd}
\]
where the vertical maps are principal $\GL_n$-bundles. In particular, $\Quot_{n,r} = Z/\GL_n$.

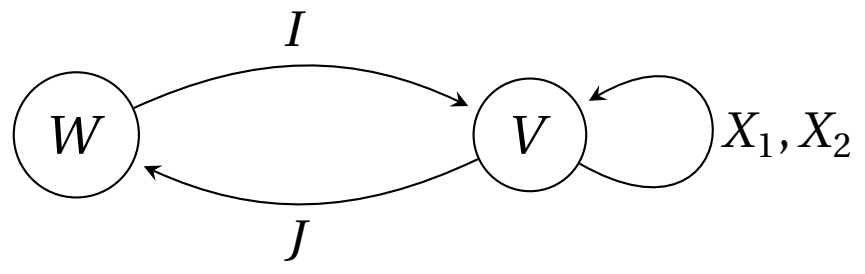
\begin{figure}[H]
\centering
\begin{tikzpicture}[>=stealth,shorten >=2pt,looseness=.5,auto]
  \matrix [matrix of math nodes,
           column sep={3cm,between origins},
           row sep={3cm,between origins},
           nodes={circle, draw, minimum size=5.5mm}]
{ 
|(A)| W & |(B)| V \\         
};
\tikzstyle{every node}=[font=\small\itshape]

\node [anchor=west,right] at (-0.25,0.7) {$I$};    \node [anchor=west,right] at (-0.25,-0.7) {$J$};  
\node [anchor=west,right] at (2.65,0) {$X_1,X_2$}; 

\draw[->] (A) to [bend left=25,looseness=1] (B) node [midway] {};
\draw[->] (B) to [bend left=25,looseness=1] (A) node [midway,below] {};
\draw[->] (B) to [out=-30,in=30,looseness=9] (B) node {};
\end{tikzpicture}
\caption{A $(r,n)$-dimensional representation of the quiver $\mathbf Q_{\mathsf{cyc}}$.}
\label{fig:cyclic-quiver}
\end{figure}

The scheme $\Quot_{n,r}$ can also be seen as a closed subscheme of the cyclic Nakajima quiver variety
\begin{equation}
\label{eqn:M^cyc}
\mathfrak{M}^{\mathsf{cyc}}_{r,n}
=\Set{(X_1,X_2,I,J)\mid [X_1,X_2]+IJ=0,\ (X_1,X_2,I)\text{ cyclic}}\big/\GL_n.
\end{equation}
Here, the starting point is the space of $(r,n)$-dimensional representations of the quiver $\mathbf Q_{\mathsf{cyc}}$ depicted in \Cref{fig:cyclic-quiver}. The cyclicity condition is identical: the image of $I$ generates the unframed representation $(X_1,X_2) \in \Rep_n(\boldit{L}_{2})$. The additional relation imposed is $[X_1,X_2] + IJ=0$. The locus in the prequotient appearing in \eqref{eqn:M^cyc} is a quasi-affine scheme which coincides with $Z$ after imposing $J=0$. After $J=0$ is imposed, the quotient gives then precisely $Z/\GL_n = \Quot_{n,r}$, which is then viewed as a closed subscheme $\Quot_{n,r} \into \mathfrak{M}^{\mathrm{cyc}}_{n,r}$.

\subsection{Symmetric products and their stratification}
\label{sec:sym-strata}
Let $Y$ be a quasiprojective variety. Denote by $\Sym^n(Y)$ its $n$-th symmetric product. There is a locally closed stratification
\begin{equation}
\label{eqn:sym-stratification}
\Sym^n(Y) = \coprod_{\lambda \vdash n} \Sym^n_\lambda(Y)
\end{equation}
indexed by partitions of $n$. We write a partition $\lambda\vdash n$ as a weakly decreasing sequence $(\lambda_1,\ldots,\lambda_l)$ of positive integers such that $\sum_i\lambda_i=n$. The \emph{length} of the partition is the number $l(\lambda) = l$. We denote by $G_{\lambda} = \{\sigma \in \FS_l \mid \lambda_{\sigma(i)}=\lambda_i \text{ for all }i\}$ the automorphism group of $\lambda$. The deepest stratum in \eqref{eqn:sym-stratification} is the \emph{small diagonal}, which is associated to the partition $\lambda = (n)$ and corresponds to the closed immersion $S \into \Sym^n(S),\,\,p\mapsto n[p].$ On the other hand, the open and dense stratum is the one corresponding to the partition $\lambda = (1,\ldots,1)$. In general, one has a locally closed immersion
\begin{equation}
\label{eqn:j-lambda-loc-closed-imm}
\map{j_\lambda}{\Sym^n_\lambda(Y)}{\Sym^n(Y),}[hook]
\end{equation}
and one can realize $\Sym^n_\lambda(Y)$ as a free quotient
\begin{align}\label{eqn: free quot Sym}
    \Sym^n_\lambda(Y) = C_\lambda(Y) / G_\lambda,
\end{align}
where $C_\lambda(Y) \subset Y^{l(\lambda)}$ is the \emph{configuration space} parametrising ordered tuples of pairwise distinct points: on this locus, the group $G_\lambda$ acts freely, and we have an induced finite \'etale cover 
\begin{equation}
\label{eqn:rholambda-finite-etale}
\map{\rho_\lambda}{C_\lambda(Y)}{\Sym^n_\lambda(Y).}
\end{equation}
Finally, we have $\Sym^{n}_{\lambda}(Y) \subseteq \overline{\Sym^{n}_{\mu}(Y)}$ if and only if $\mu \leq \lambda$, with respect to the natural partial order on partitions.

\begin{lemma}
\label{lemma:addition-is-etale} 
Fix an integer $k\geq 2$ and a variety $Y$ with $\dim Y>0$. Let $\boldit{n} = (n_1,\ldots,n_k)$ be a $k$-tuple of nonnegative integers such that $n_1+\cdots+n_k = n$. Consider the addition of cycles map
\[
\Map{\cup_{\boldit n}}{\displaystyle\prod_{1\leq i\leq k}\Sym^{n_i}(Y)}{\Sym^n(Y).}
\]
Let $A_{\boldit n}$ be the open subset of the source parametrising $k$-tupes $(\xi_1,\ldots,\xi_k)$ of cycles such that $\Supp \xi_i \cap \Supp \xi_j = \emptyset$ for all $i \neq j$. Then the restricted summation map
\[
\Map{\cup_{\boldit n}}{A_{\boldit n}}{\Sym^n(Y)}
\]
is \'etale. For $k=2$, the images of these maps cover the complement $\Sym^n(Y) \setminus Y$ of the small diagonal.
\end{lemma}

\begin{proof}
Let $\CA_{\boldit n} \subset \prod_{1\leq i\leq k} \CCoh^{n_i}(Y)$ be the open substack of $k$-tuples of sheaves with pairwise disjoint supports. Then, by \cite[Cor.~1.35]{Fantechi-Ricolfi-structural}, the direct sum map $\CA_{\boldit n} \to \CCoh^n(Y)$ is \'etale. We have a commutative diagram
\[
\begin{tikzcd}[column sep=large]
\CA_{\boldit n}\arrow{d}\arrow{r}{\oplus} & \CCoh^n(Y)\arrow{d}{\sigma_n} \\
A_{\boldit n}\arrow{r}{\cup_{\boldit n}} & \Sym^n(Y)
\end{tikzcd}
\]
where the vertical morphisms are good moduli space morphisms. We need to prove that the bottom map is \'etale. Since both source and target are finite type $\BC$-schemes, it is enough to check that the induced map on completed local rings at closed points is an isomorphism. Say $\xi = \cup_{\boldit n}(\xi_\bullet)$ for a $k$-tuple $\xi_\bullet = (\xi_1,\ldots,\xi_k)\in A_{\boldit n}$. Let $E_\bullet = (E_1,\ldots,E_k) \in \CA_{\boldit n}$ be the closed point above $\xi_\bullet$, and let $E \in \CCoh^n(Y)$ be its image along $\oplus$, so that $E$ lies above the cycle $\xi$. By definition of $\CA_{\boldit n}$, the map $\oplus$ is stabilizer-preserving, in particular the natural map $\prod_{i}\Aut(E_i) \to \Aut(E)$ is an isomorphism. Moreover, the natural map $\mathsf A_E \to \mathsf A_{E_\bullet}$ between formal miniversal rings is an isomorphism by \'etaleness. Taking stabilizer-invariants then yields an isomorphism
\[
\Map{\sim}{\mathsf A_E^{\Aut(E)}}{\displaystyle\widehat{\bigotimes_{1\leq i\leq k}} \mathsf A_{E_i}^{\Aut(E_i)}.}
\]
By the \'etale-local structure of good moduli spaces, cf.~\cite[Lemma 2.7]{Chen-Pertusi-Zhao} and \cite{Alper-local-properties}, this isomorphism is precisely the map $\widehat{\OO}_{\Sym^n(Y),\xi} \to \widehat{\OO}_{A_{\boldit n},\xi_\bullet}$ we were interested in.
\end{proof}

\subsection{\texorpdfstring{Perverse sheaves on $\Sym^n(S)$}{}}
Let $S$ be a smooth quasiprojective surface and $\SE$ a  locally free sheaf on $S$. Given the Quot-to-Chow morphism 
\[
\Map{\pi_{n,\SE}}{\Quot^n(\SE) }{\Sym^n(S),}
\]
define the mixed Hodge module
\begin{align}\label{eqn: def Fn}
\CF_{n,\SE} \defeq \pi_{n,\SE\ast}\mathbb{Q}_{\Quot^n(\SE)} \,\in\,\Db(\MHM(\Sym^n(S))).    
\end{align}

Let us denote by $\CT$ the stratification \eqref{eqn:sym-stratification}.

\begin{proposition} \label{sheafisconstructible}
Let $S$ be a smooth quasiprojective surface and $\SE$ a locally free sheaf on $S$ of rank $r\geq 1$. Then $\mathcal{F}_{n,\SE}$ is $\CT$-constructible, i.e. 
\[ 
\mathcal{F}_{n,\SE} \in \Db_{\CT}(\MHM(\Sym^n(S))).
\]
\end{proposition}

\begin{proof}
Fix a partition $\lambda \vdash n$.
Consider the morphisms
\[
\begin{tikzcd}[row sep=large]
C_\lambda(S) \arrow{d}{\rho_\lambda} & \\
\Sym^n_\lambda(S) \arrow[hook]{r}{j_\lambda} & \Sym^n(S)
\end{tikzcd}
\]
introduced in \Cref{sec:sym-strata}, with $j_\lambda$ locally closed and $\rho_\lambda$ finite \'etale.
The fibre of $\pi_{n,\SE}$ over a cycle $\xi = \sum_{i}\lambda_i[p_i] \in \Sym^n_\lambda(S)$ is
\[
\Quot_{n,r}^\lambda \defeq \prod_{i=1}^{l(\lambda)} \Quot_{\lambda_i,r}^0,
\]
where the factors are punctual Quot schemes, cf.~\Cref{def:punctual-quot}.
We have a cartesian diagram 
\[
\begin{tikzcd}
D^{\lambda}_{n,\mathcal{E}}\defeq C_\lambda(S) \times_{\Sym^n(S)} \Quot^n(\SE) \MySymb{dr}\arrow[r] \arrow[d,"q"']
& \Quot^n(\SE)  \arrow[d,"\pi_{n,\SE}"] \\
C_\lambda(S) \arrow[r,"j_\lambda\rho_\lambda"'] & \Sym^n(S)
\end{tikzcd},
\]
where the base change map $q$ is a Zariski trivial fibration with fibres isomorphic to $\Quot_{n,r}^\lambda$ (see e.g.~the proof of \cite[Thm.~2.3]{Ricolfi-Power-structure-Quot}). By étale base change,
\[ 
\rho_{\lambda}^*(j_{\lambda}^*\mathcal{F}_{n,\SE}) \simeq q_{*}\mathbb{Q}_{D^{\lambda}_{n,\mathcal{E}}}.
\]
In particular,  $\rho_{\lambda}^*(j_{\lambda}^*\mathcal{F}_{n,\SE})$ is a complex whose cohomology sheaves are local systems. Finally, since $\rho_{\lambda}$ is surjective and étale, it follows that $j_{\lambda}^*\mathcal{F}_{n,\SE}$ is also a complex whose cohomology sheaves are local systems, which concludes the proof.
\end{proof}

\subsection{Factorization structure}\label{subsec: factorization str}

Let $S$ be a smooth surface. For $k \geq 2$ and $\boldit n = (n_1,\ldots,n_k)$, recall from \Cref{lemma:addition-is-etale} the open subscheme 
\[
\map{\varepsilon_{\boldit n}}{A_{\boldit n}}{\Sym^{\boldit n}(S)\defeq\displaystyle\prod_{i=1}^k \Sym^{n_i}(S)}[hook]
\]
parametrising $k$-tuples of cycles with pairwise disjoint support. Recall that the addition of cycles $\cup_{\boldit n} \colon A_{\boldit n} \to \Sym^n(S)$ is \'etale. In this context, we can prove the following result.

\begin{proposition}
\label{prop:sheaffactorization}
 For a locally free sheaf $\SE$ on $S$, there is an isomorphism of complexes of mixed Hodge modules
\[
\cup_{\boldit n}^\ast \CF_{n,\SE} \cong \varepsilon_{\boldit n}^\ast (\CF_{n_1,\SE} \boxtimes \cdots \boxtimes \CF_{n_k,\SE})\,\in\,\Db(\MHM(A_{\boldit n})).
\]
\end{proposition}

\begin{proof}
Set $\Quot^{\boldit{n}}(\SE) \defeq \prod_i\Quot^{n_i}(\SE)$ and denote by $\pi_{\boldit n,\SE} \colon \Quot^{\boldit{n}}(\SE) \to \Sym^{\boldit n}(S)$ the product $\pi_{n_1,\SE} \times \cdots \times \pi_{n_k,\SE}$ of Quot-to-Chow maps.
One observes that both squares in the commutative diagram
\[
\begin{tikzcd}[row sep=large,column sep=large]
\Quot^{\boldit{n}}(\SE) & \pi_{\boldit n,\SE}^{-1}(A_{\boldit n}) & \Quot^n(\SE)  \\
\Sym^{\boldit n}(S) & A_{\boldit n} & \Sym^n(S)
\arrow["\pi_{\boldit n,\SE}", from=1-1, to=2-1]
\arrow[from=1-2, to=1-1,hook']
\arrow[from=1-2, to=1-3]
\arrow[swap,"\rho"', from=1-2, to=2-2]
\arrow["{\pi_{n,\SE}}", from=1-3, to=2-3]
\arrow[swap,"{\varepsilon_{\boldit n}}", from=2-2, to=2-1,hook']
\arrow[swap,"{\cup_{\boldit n}}"', from=2-2, to=2-3]
\end{tikzcd}
\] 
are cartesian.
This is because a torsion sheaf with support on $\sum_i\xi_i$ where $\Supp \xi_i \cap \Supp \xi_j=\emptyset$ for all $i\neq j$ splits into $k$ torsion sheaves supported on cycles $\xi_1,\ldots,\xi_k$. By base change applied to both squares, we have
\[
\cup_{\boldit n}^\ast \CF_{n,r}
\cong \rho_\ast \BQ_{\pi_{\boldit n,\SE}^{-1}(A_{\boldit n})} 
\cong \varepsilon_{\boldit n}^\ast \pi_{\boldit n,\SE\ast}\BQ_{\Quot^{\boldit{n}}(\SE)} \cong \varepsilon_{\boldit n}^\ast (\CF_{n_1,\SE} \boxtimes \cdots \boxtimes \CF_{n_k,\SE}),
\]
as required.
\end{proof}

\section{Cohomological purity}\label{sec: hyperbloc}

\subsection{Hyperbolic localization}\label{subsec: Hyp Loc}
The closed immersion $\BA^1 \times \set{0} \into \BA^2$ induces a closed immersion 
\[
\begin{tikzcd}
\mathrm{SN}^n(\BA^2) \defeq  \Sym^n(\BA^1 \times \set{0}) \arrow[hook]{r} & \Sn,
\end{tikzcd}
\]
where the source, which we refer to as the \emph{semi-nilpotent} part of the symmetric product, is the locus consisting of cycles in $\Sn$ supported on $\BA^1 \times \set{0}$. We denote by
\[
\begin{tikzcd}
\Quot^{\mathrm{SN}}_{n,r}=\pi_{n,r}^{-1}(\mathrm{SN}^n(\BA^2))\arrow[hook]{r} & \Quot_{n,r}
\end{tikzcd}
\]
its preimage along the Quot-to-Chow morphism $\pi_{n,r}$. 
Consider the $\BG_m$-action on $\BA^2$ given by 
\[
t\cdot(x,y)\defeq  (tx,t^{-1}y).
\]
This torus action induces a $\BG_m$-action on $\Sn $, whose only fixed point is $n[0]\in \Sn$ and whose attracting locus is $\mathrm{SN}^n(\BA^2)$.
In other words, we have maps
\[ 
\begin{tikzcd}
n[0] \arrow[hook]{r}{f_n} & \mathrm{SN}^n(\BA^2) \arrow[hook]{r}{g_n} & \Sn,    
\end{tikzcd}
\]
where $f_n$ is the inclusion of the fixed locus and $g_n$ is  the inclusion of the attracting locus.

Define the hyperbolic localization functor 
\[ 
\begin{tikzcd}
\Db_{\CT}(\MHM(\Sym^n(\BA^2)))  \arrow{r}{\Phi_n} & \Db(\mathrm{MHS}), \quad \mathcal{F} \mapsto f_n^{!}g_n^*\mathcal{F}.
\end{tikzcd}
\]
The functor $\Phi_n$ is a lift of the \emph{hyperbolic localization} functor considered in \cite{bradenmautner}.

\begin{proposition} 
\label{preservetstructure}
The functor $\Phi_n$ preserves t-structures. \end{proposition}
\begin{proof}
Since $\Phi_n$ is defined using the six-functor formalism, it commutes with the rationalization functor $\mathrm{rat}_{\Sn}$. It is proved in \cite[Prop.~4.3]{bradenmautner} that the hyperbolic localization functor 
\[
\begin{tikzcd}
\mathrm{rat}_{\Sn}(\Phi_n)\colon  \Db_{\CT}(\Sn)  \arrow{r} &  \Db(\mathrm{Vec}_{\mathbb{Q}})
\end{tikzcd}
\] preserves t-structures. Furthermore, the rationalization functor is conservative and hence the claim follows for mixed Hodge modules as well. 
\end{proof}

The next result will allow us to check the vanishing required for purity \emph{after} the application of the functor $\Phi_n$.

\begin{proposition} 
\label{onthediagonal}
    Let $\mathcal{G} \in \mathrm{MHM}_{\CT}(\Sym^n(\BA^2))$ be a mixed Hodge module supported on the small diagonal $\BA^2\into \Sym^n(\BA^2)$ such that $\Phi_n(\mathcal{G})=0$. Then $\mathcal{G}=0$. 
\end{proposition}
\begin{proof}
Let $i_n\colon  \BA^2 \hookrightarrow \Sn$ denote the  closed immersion of the small diagonal. Set $\mathcal{G}= i_{n\ast}\mathcal{G}^{\prime}$, where $\mathcal{G}^{\prime}$ is a mixed Hodge module on $\BA^2$. By definition of the stratification $\CT$, we have that  $\mathrm{rat}_{\BA^2}(\mathcal{G}^{\prime})$ is a constant local system, since $\BA^2$ is simply connected. Therefore, we have that $\mathcal{G}^{\prime} =
H[2]$, where $H$ is the pullback of a constant mixed Hodge structure  along $\BA^2\to \pt$.

We conclude by noticing that
 \[ 
 \Phi_n(\mathcal{G}) \simeq \Phi_n(i_{n\ast}H[2]) \simeq \HH_c(\BA^1,\mathbb{Q})[2] \otimes H,
 \]
 which implies  $\Phi_n(\mathcal{G})=0$ if and only if $\mathcal{G}=0$.
\end{proof}
\subsection{\texorpdfstring{Purity of $\CF_{n,\SE}$}{}}
\label{subsec: proof of main thm}
The goal of this subsection is to prove \Cref{thm:Main thm inro} from the introduction.
\begin{proposition} \label{purityondiagonal}
Fix integers $r\geq 1$ and $n\geq 0$. Then the mixed Hodge structure on   the compactly supported cohomology 
\[
\mathrm{H}_c^k(\Quot^{\mathrm{SN}}_{n,r})
\] is pure of weight $k$.
\end{proposition} 
\begin{proof}
We have a cartesian diagram
\[
\begin{tikzcd}
\Quot^{\mathrm{SN}}_{n,r} \MySymb{dr} \arrow[hook]{r}\arrow{d} & \Quot_{n,r} \arrow{d}{\pi_{n,r}} \\
\mathrm{SN}^n(\BA^2)\arrow[hook]{r} & \Sym^n(\BA^2)
\end{tikzcd}
\]
defining the semi-nilpotent Quot scheme. Consider the dual quiver $\mathbf Q_{\mathsf{cyc}}^\vee$ to the one considered in \Cref{fig:cyclic-quiver}, depicted in \Cref{fig:cyclic-quiver-dualized}.
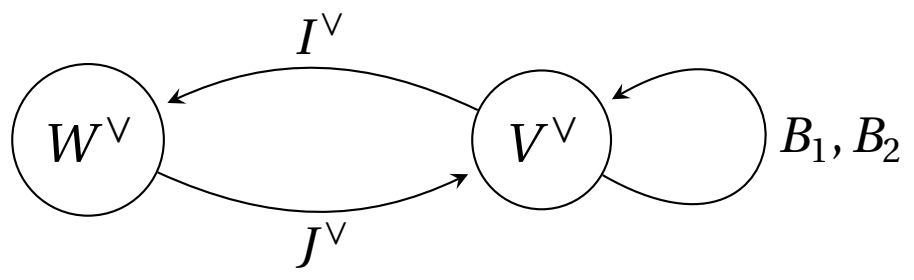
\begin{figure}[H]
\centering
\begin{tikzpicture}[>=stealth,shorten >=2pt,looseness=.5,auto]
  \matrix [matrix of math nodes,
           column sep={3cm,between origins},
           row sep={3cm,between origins},
           nodes={circle, draw, minimum size=5.5mm}]
{ 
|(A)| W^\vee & |(B)| V^\vee \\         
};
\tikzstyle{every node}=[font=\small\itshape]

\node [anchor=west,right] at (-0.25,0.7) {$I^\vee$};    \node [anchor=west,right] at (-0.25,-0.7) {$J^\vee$};  
\node [anchor=west,right] at (2.95,0) {$B_1,B_2$}; 

\draw[->] (B) to [bend right=25,looseness=1] (A) node [midway] {};
\draw[->] (A) to [bend right=25,looseness=1] (B) node [midway,below] {};
\draw[->] (B) to [out=-30,in=30,looseness=9] (B) node {};
\end{tikzpicture}
\caption{A $(r,n)$-dimensional representation of the dual quiver $\mathbf Q_{\mathsf{cyc}}^\vee$.}
\label{fig:cyclic-quiver-dualized}
\end{figure}
Consider the quasi-affine scheme
\[
Z' \defeq \Set{(B_1,B_2,I^\vee)|[B_1,B_2]=0,\,B_1\text{ is nilpotent,}\,(B_1,B_2,I^\vee)\text{ is cocyclic}}\subset  \Rep_{(r,n)}(\mathbf Q_{\mathsf{cyc}}^\vee).
\]
Here, \emph{cocyclic} means that no nontrivial subspace of $V^\vee$ preserved by $B_1$ and $B_2$ is contained in $\ker I^\vee$. We have a natural $\GL_n$-equivariant morphism 
\[
\Map{\dd}{Z'}{Z}, \qquad (B_1,B_2,I^\vee)\mapsto (B_2^\vee,B_1^\vee,(I^\vee)^\vee)
\]
where 
\[
Z=\Set{(X_1,X_2,I)|[X_1,X_2]=0,\,(X_1,X_2,I)\text{ is cyclic}}\subset  \Rep_{(r,n)}(\mathbf Q_{\mathsf{cyc}})
\]
is the quasi-affine scheme defined in \Cref{sec:quot-via-quivers}. Consider now the \emph{semi-nilpotent} quiver variety $\mathfrak{L}^0_{n,r}$ associated to the framed Jordan quiver, as defined (more generally) in \cite[Def.~4.1.1]{schiffmannvasserot}. Explicitly, in our case we have $\mathfrak{L}^0_{n,r} \defeq Z'/\GL_n$.
The morphism $\dd$ then descends to a morphism 
\[
\map{\overline{\dd}}{\mathfrak L^0_{n,r}}{Z/\GL_n = \Quot_{n,r}}
\]
whose composition with $\pi_{n,r}$ factors through $\mathrm{SN}^n(\BA^2) \into \Sym^n(\BA^2)$ by the nilpotency of $B_1$. The canonical dualising morphism induced by the universal property of fibre products
\[
\map{\widetilde\dd}{\mathfrak L^0_{n,r}}{\Quot^{\mathrm{SN}}_{n,r}}
\]
is easily checked to be an isomorphism at the level of reduced underlying schemes.

Our conclusion then follows from \cite[Prop.~4.3]{schiffmannvasserot}, which proves that the compactly supported cohomology  $\HH^k_c(\mathfrak{L}^{0}_{n,r})$ is pure of weight $k$. 
\end{proof}

We are now ready to prove \Cref{thm:Main thm inro} from the Introduction.
\begin{theorem}
\label{main-thm-purity-A2}
Let $S$ be a smooth quasiprojective surface and $\SE$ a locally free sheaf on $S$.   Then the mixed Hodge module 
\[\CF_{n,\SE} \defeq \pi_{n,\SE\ast}\mathbb{Q}_{\Quot^n(\SE)} \,\in\,\Db(\MHM(\Sym^n(S))) \] 
is pure of weight $0$.
\end{theorem}
\begin{proof}
Set  $r=\rk \SE$. By \Cref{lemma:purity-checked-on-etale-cover}, purity of $\CF_{n, \SE}$ can be checked étale-locally, therefore we may assume that $(S, \SE)=(\BA^2, \OO^{\oplus r}_{\BA^2})$.
By \Cref{preservetstructure}, we have that
\begin{align}\label{eqn: Hk}
\Phi_n({}^{\mathfrak{p}}\mathcal{H}^k(\mathcal{F}_{n,r})) \simeq \mathrm{H}^k(\Phi_n(\mathcal{F}_{n,r})), \qquad \mbox{for all } k\in\mathbb{Z}.
\end{align}
Let us also notice that $\mathcal{F}_{n,r}$ is naturally a $\BG_m$-equivariant complex. Therefore, by \cite[Thm.~8]{bradenhyperbolic}, the functor $\Phi_n$ preserves the purity and weight of the $\BG_m$-equivariant pure mixed Hodge modules $\Gr_w({}^{\mathfrak{p}}\mathcal{H}^k(\mathcal{F}_{n,r}))$, i.e.
\begin{align}\label{eqn: Gr}
   \Phi_n(\Gr_w({}^{\mathfrak{p}}\mathcal{H}^k(\mathcal{F}_{n,r}))) \simeq \Gr_w (\Phi_n({}^{\mathfrak{p}}\mathcal{H}^k(\mathcal{F}_{n,r}))), \qquad \mbox{for all } w\in\mathbb{Z}.
\end{align}
 Denote by $p_{n}\colon  \mathrm{SN}^n(\BA^2) \rightarrow n[0]$ the projection. By \cite[Thm.~2.10.3]{achar} we have the isomorphism
\[
f_n^!g_n^*\mathcal{F}_{n,r} \simeq p_{n!}g_n^*\mathcal{F}_{n,r},
\] 
where we exploited once more the $\BG_m$-equivariance of $\CF_{n,r}$. In particular, we have an identification
\begin{equation}
\label{eqn: SN}
\mathrm{H}^k(\Phi_n(\mathcal{F}_{n,r}) )= \mathrm{H}_c^k(\Quot^{\mathrm{SN}}_{n,r}), \qquad \mbox{for all } k\in\mathbb{Z}.
\end{equation} 
Combining \eqref{eqn: Hk}, \eqref{eqn: Gr} and \eqref{eqn: SN} we conclude that
\begin{equation}\label{eq:1st way of computing hyperbolic localization}
    \Phi_n(\Grh) \simeq \Gr_w(\mathrm{H}_c^k(\Quot^{\mathrm{SN}}_{n,r})), \qquad \mbox{for all } w,k\in\mathbb{Z}.
\end{equation}
By \Cref{purityondiagonal} the mixed Hodge structure on $\mathrm{H}_c^k(\Quot^{\mathrm{SN}}_{n,r})$ is pure of weight $k$, which combined with \eqref{eq:1st way of computing hyperbolic localization} yields the  vanishing 
\begin{align}\label{eqn: final zero}
    \Phi_n(\Grh)=0, \qquad \mbox{for all } w\neq k\in\mathbb{Z}.
\end{align}
We  claim now that
\[
\Grh=0, \qquad \mbox{for all } w\neq k\in\mathbb{Z},
\]
by which we will conclude purity of $\CF_{n,r}$. We prove the claimed vanishing by induction on $n\geq 1$. 

For $n=1$, the Quot-to-Chow morphism is simply
\begin{align*}
\Map{\pi_{1,r}}{\Quot_{1,r}\cong \BP^{r-1}\times \BA^2}{\BA^2\cong \Sym^1(\BA^2)},
\end{align*}
which implies that the mixed Hodge module 
\begin{align*}
    \CF_{1,r}=\HH^*(\BP^{r-1})
\end{align*}
is pure of weight $0$. 

Fix now $n\geq 2$.  \Cref{lemma:addition-is-etale}, applied with $k=2$, implies that the morphisms $\cup_{\boldit n}\colon A_{\boldit n} \to \Sym^n(\BA^2)$, ranging over pairs $\boldit n = (n_1,n_2)$ with $n_1+n_2=n$, provide an étale cover of the complement $\Sym^n(\BA^2) \setminus \BA^2$ of the small diagonal. Thus, by \Cref{lemma:purity-checked-on-etale-cover}, the vanishing 
\begin{align}\label{eqa vanishing 1}
    \Grh|_{\Sym^n(\BA^2) \setminus \BA^2}=0, \qquad \mbox{for all } w\neq k\in\mathbb{Z},
\end{align}
is equivalent to the vanishing 
\begin{align}\label{eq: vanishing 2}
    \cup_{\boldit n}^\ast \Grh=0,  \qquad \mbox{for all } w\neq k\in\mathbb{Z}\,\,\, \mbox{and for all } \boldit{n}=(n_1, n_2).
\end{align}
The morphism $\cup_{\boldit n}$ is étale, therefore the pullback $\cup_{\boldit n}^*$ preserves t-structures and the weight filtration. Therefore, \Cref{prop:sheaffactorization} yields
\begin{align*}
    \cup_{\boldit n}^* \Grh &\simeq \Gr_{w}({}^{\mathfrak{p}}\mathcal{H}^k(\cup_{\boldit n}^* \mathcal{F}_{n,r}))\\
    &\simeq \Gr_{w}({}^{\mathfrak{p}}\mathcal{H}^k(\varepsilon_{\boldit n}^*(\mathcal{F}_{n_1,r} \boxtimes \mathcal{F}_{n_2,r})))\\
    &\simeq \varepsilon_{\boldit n}^* \Gr_{w}({}^{\mathfrak{p}}\mathcal{H}^k(\mathcal{F}_{n_1,r} \boxtimes \mathcal{F}_{n_2,r}))\\
    &\simeq \varepsilon_{\boldit n}^* \bigoplus_{\substack{i+j=k\\a+b=w}}\Gr_{a}({}^{\mathfrak{p}}\mathcal{H}^i(\mathcal{F}_{n_1,r} ))\boxtimes \Gr_{b}({}^{\mathfrak{p}}\mathcal{H}^j(\mathcal{F}_{n_2,r}))
\end{align*}
where in the third line we used that $ \varepsilon_{\boldit n}$ is an open embedding. Clearly,  if $w\neq k$, then $(i,j)\neq (a,b)$ for each summand.  Assume without loss of generality that $ i\neq a$. Since $n_1<n$, by the induction hypothesis we have the vanishing
\begin{align*}
    \Gr_{a}({}^{\mathfrak{p}}\mathcal{H}^i(\mathcal{F}_{n_1,r} ))=0,
\end{align*}
which implies the vanishing \eqref{eq: vanishing 2} for all $\boldit{n}=(n_1, n_2)$ and in turn the vanishing \eqref{eqa vanishing 1}. In particular, the mixed Hodge module $\Grh$ is supported on the small diagonal $\BA^2\into \Sym^n(\BA^2)$. We conclude the induction step by combining the vanishing \eqref{eqn: final zero} with \Cref{onthediagonal}.
\end{proof}

\section{Semisimple decomposition}\label{section:computation of decomposition}
Let $S$ be a smooth quasiprojective surface and $\SE$ a locally free sheaf on $S$. By \Cref{main-thm-purity-A2}, the complex of mixed Hodge modules 
\[
\mathcal{F}_{n,\SE} = \pi_{n,\SE\ast}\mathbb{Q}_{\Quot^n(\SE)} \,\in\, \Db_{\CT}(\MHM(\Sym^n(S)))
\] 
is pure of weight $0$. Since the category of pure mixed Hodge modules is semisimple \cite[Cor.~5.2.13]{Saitofrench},
by \cite[Thm.~3.21]{saito} there exists a (non-canonical) isomorphism
\begin{equation}
\label{splitting0}
\mathcal{F}_{n,\SE} \simeq \bigoplus_{\substack{\lambda \vdash n \\ j \in \mathbb{Z}}} \mathrm{IC}_{\overline{\Sym^n_{\lambda}(S)}}(\mathbb{V}_{\lambda,j})[-j],  \end{equation}
where   $\mathbb{V}_{\lambda,j}$ is a variation of mixed Hodge structures, pure of weight $j-2l(\lambda)$, on $\Sym^n_{\lambda}(S)$. Since $\Sym^n_{\lambda}(S)$ is smooth of dimension $2l(\lambda)$, we have that
\begin{equation}
\label{IC-restricted-lambda}
\mathrm{IC}_{\overline{\Sym^n_{\lambda}(S)}}(\mathbb{V}_{\lambda,j})\big|_{\Sym^n_{\lambda}(S)} \simeq \mathbb{V}_{\lambda,j}[2 l(\lambda)].
\end{equation}

\begin{remark}
Göttsche--Soergel's classical result \cite{gotsoergelhilbert} computes the decomposition \eqref{splitting0} in the case $\SE=\OO_S$, that is for the Hilbert scheme of points  $\Hilb^n(S)$. In the rank 1 case, the Hilbert scheme is smooth of dimension $2n$ and the Quot-to-Chow morphism $\pi_{n, \OO_S}$ is a proper semismall map. In this setting, the decomposition theorem applies and shows that $\CF_{n, \OO_S}=\pi_{n, \OO_S\ast}\BQ_{\Hilb^n(S)}$ is a shift of a perverse sheaf and, in particular, each direct summand corresponding to $j\neq 2n$ vanishes. In our case, that is for a general higher rank locally free sheaf $\SE$, the Quot scheme $\Quot^n(\SE)$ is not smooth and the Quot-to-Chow map $\pi_{n, \SE}$ is not semismall, thus the argument in \cite{gotsoergelhilbert} cannot be  generalized.
\end{remark}
Recall from \eqref{eqn: free quot Sym} that $\Sym^n_\lambda(S) = C_\lambda(S) / G_\lambda$, and denote by 
\[ 
\Map{\rho_{\lambda}}{C_{\lambda}(S)}{\Sym^n_{\lambda}(S)} 
\]
the associated finite étale morphism \eqref{eqn:rholambda-finite-etale}.
 Set
\begin{align*}
     \mathbb{H}_{\lambda,j} \defeq  \rho_{\lambda}^* \mathbb{V}_{\lambda,j}[2l(\lambda)],
\end{align*}
which is naturally a $G_\lambda$-equivariant complex.
\begin{proposition} \label{prop:etalepullbackofV}
Let $\lambda$ be a partition of $n$ and $j\in \BZ$. 
There is an isomorphism 
\[
(\rho_{\lambda\,\ast} \mathbb{H}_{\lambda,j})^{G_{\lambda}}  \simeq \mathbb{V}_{\lambda,j}[2l(\lambda)].
\]
\end{proposition}
\begin{proof}
The proof follows directly from the fact that finite group $G_\lambda$ acts freely on $  C_\lambda(S) $, which implies that the category of $G_\lambda$-equivariant variations of mixed Hodge structures on $  C_\lambda(S)$ is equivalent to the category of  variations of mixed Hodge structures on $\Sym_{\lambda}^n(S) $, see e.g.~\cite{Bernstein-Lunts}.
\end{proof}

We show now that the mixed Hodge module $\mathbb{H}_{\lambda,j}$ is actually the pullback of a constant mixed Hodge structure.
\begin{proposition}\label{Prop:MHS on Cl are constant}
Let $\lambda$ be a partition of $n$ and $j\in \BZ$. Then $\mathbb{H}_{\lambda,j}$ is constant, i.e.~we have an isomorphism
    \[ \mathbb{H}_{\lambda,j} \simeq H_{\lambda,j},
    \]
    where $H_{\lambda,j}$ is pullback of a constant mixed Hodge structure $H_{\lambda,j}$ to $C_{\lambda}(S)$. 
\end{proposition}

\begin{proof}
Pulling back along $j_\lambda \colon \Sym^n_\lambda(S) \into \Sym^n(S)$ the decomposition in \eqref{splitting0} and using \eqref{IC-restricted-lambda} yields
\begin{align} 
\label{splitting}
j_{\lambda}^* \mathcal{F}_{n,\SE} &\simeq \bigoplus_{\substack{\mu \vdash n \\ j \in \mathbb{Z}}} j_{\lambda}^* \mathrm{IC}_{\overline{\Sym^n_{\mu}(S)}}(\mathbb{V}_{\mu,j})[-j]  \\ & \simeq \bigoplus_{j\in \mathbb{Z}}\mathbb{V}_{\lambda,j}[2l(\lambda)-j] \oplus \bigoplus_{\substack{\mu\neq \lambda \\ j \in \mathbb{Z}}} j_{\lambda}^* \mathrm{IC}_{\overline{\Sym^n_{\mu}(S)}}(\mathbb{V}_{\mu,j})[-j]\nonumber,
\end{align}
 which implies that $\mathbb{H}_{\lambda,j}[-j]=\rho_{\lambda}^*\mathbb{V}_{\lambda,j}[2l(\lambda)-j]$ is a direct summand of $\rho_{\lambda}^*j_{\lambda}^*\mathcal{F}_{n,\SE}$. Consider the cartesian diagram
 \[
\begin{tikzcd}[row sep=large]
	{D^{\lambda}_{n,\mathcal{E}}} && {\Quot^n(\SE)} \\
	{C_{\lambda}(S)} & {\Sym^n_{\lambda}(S)} & {\Sym^n(S)}
	\arrow[from=1-1, to=1-3]
	\arrow["{\pi_{\lambda}}"', from=1-1, to=2-1]
	\arrow["\pi_{n,\SE}",from=1-3, to=2-3]
	\arrow["{\rho_{\lambda}}", from=2-1, to=2-2]
	\arrow["{j_{\lambda}}",hook, from=2-2, to=2-3]
\end{tikzcd}
\]
defining the scheme $D^{\lambda}_{n,\mathcal{E}}$.
By base change we have that $\rho_{\lambda}^*j_{\lambda}^*\mathcal{F}_{n,\SE} \simeq \pi_{\lambda\ast}\mathbb{Q}_{D^{\lambda}_{n,\mathcal{E}}}$.
Since the morphism $\pi_{\lambda}$ is a Zariski locally trivial fibration (see the proof of \cite[Thm.~2.3]{Ricolfi-Power-structure-Quot}), the cohomology sheaves of the pushforward $\pi_{\lambda\ast}\mathbb{Q}_{D^{\lambda}_{n,\mathcal{E}}}$ are constant variations of mixed Hodge structures on $C_{\lambda}(S)$. Since $\mathbb{H}_{\lambda,j}[-j]$ is concentrated in a single degree and is a direct summand of $\pi_{\lambda\ast}\mathbb{Q}_{D^{\lambda}_{n,\mathcal{E}}}$, we conclude that $\mathbb{H}_{\lambda,j}$ is a constant variation of a mixed Hodge structures.  
\end{proof}

Let us now consider the diagram
\begin{equation}\label{Diagram:Diagram with maps}
\begin{tikzcd}[row sep=large,column sep=large]
	{S^{l(\lambda)}} & {\overline{\Sym^n_{\lambda}(S)}} & {\Sym^n(S)} \\
	{C_{\lambda}(S)} & {\Sym^n_{\lambda}(S)}
	\arrow[from=1-1, to=1-2]
	\arrow[swap,"{q_{\lambda}}"', bend left=30, from=1-1, to=1-3]
	\arrow[hook,from=1-2, to=1-3]
	\arrow["{\eta_{\lambda}}", hook, from=2-1, to=1-1]
	\arrow["{\rho_{\lambda}}", from=2-1, to=2-2]
	\arrow["{i_{\lambda}}"',hook, from=2-2, to=1-2]
	\arrow["{j_{\lambda}}"', hook, from=2-2, to=1-3]
\end{tikzcd}
\end{equation}
where we define the morphism $q_{\lambda}:S^{l(\lambda)}\rightarrow \Sym^n(S)$ by
\[
(p_1,\ldots,p_l)\mapsto \sum_{i=1}^l \lambda_i [p_i].
\]
Recall that, by \Cref{Prop:MHS on Cl are constant}, we have that
   \[ \mathbb{H}_{\lambda,j} \simeq H_{\lambda,j}
    \]
    is a constant mixed Hodge structure.
\begin{proposition}\label{Prop:ICs as equi sheaves}
Let $\lambda$ be a partition of $n$ and $j\in \BZ$. We have an isomorphism 
    \[ 
    \left(q_{\lambda\ast} H_{\lambda,j}\right)^{G_{\lambda}} \simeq \mathrm{IC}_{\overline{\Sym^n_{\lambda}(S)}}(\mathbb{V}_{\lambda,j}).
    \] 
    In particular, we have
    \[ 
    \mathcal{F}_{n,\SE} \simeq \bigoplus_{\substack{\lambda \vdash n \\ j \in \mathbb{Z}}} \,\left(q_{\lambda\ast} H_{\lambda,j}\right)^{G_{\lambda}}[-j].
    \]
\end{proposition}

\begin{proof}
Since $S^{l(\lambda)}$ is smooth and $C_{\lambda}(S)$ is dense open in $S^{l(\lambda)}$, we have that $(\eta_{\lambda})_{!*}(H_{\lambda,j}) \simeq H_{\lambda,j}$ (see \cite[Section 1.13]{saitonotes}).
Since the morphism $q_{\lambda}$ is finite, the functor $q_{\lambda\ast}$ is t-exact and thus
\begin{align*}
q_{\lambda\ast} H_{\lambda,j} \simeq q_{\lambda\ast} \eta_{\lambda!*}(H_{\lambda,j}) &\simeq q_{\lambda\ast} \mathrm{Im}(\pHcal^0(\eta_{\lambda!} H_{\lambda,j}) \rightarrow \pHcal^0(\eta_{\lambda \ast}H_{\lambda,j}))   \\ & \simeq \mathrm{Im}(\pHcal^0((q_{\lambda}\eta_{\lambda})_! H_{\lambda,j}) \rightarrow \pHcal^0((q_{\lambda}\eta_{\lambda})_*H_{\lambda,j})) \\ & \simeq  \mathrm{Im}(\pHcal^0((j_{\lambda}\rho_{\lambda})_! H_{\lambda,j}) \rightarrow \pHcal^0((j_{\lambda}\rho_{\lambda})_{*}H_{\lambda,j})) \\ &\simeq j_{\lambda!*}(\rho_{\lambda\ast} H_{\lambda,j}).
\end{align*}
Taking $G_\lambda$-invariants yields
\[
\left(q_{\lambda\ast}H_{\lambda,j}\right)^{G_{\lambda}} \simeq \left(j_{\lambda!*}\rho_{\lambda\ast}H_{\lambda,j}\right)^{G_{\lambda}} \simeq j_{\lambda!*}\left(\rho_{\lambda\ast}H_{\lambda,j}\right)^{G_\lambda}=\mathrm{IC}_{\overline{\Sym^n_{\lambda}(S)}}(\mathbb{V}_{\lambda,j}),
    \]
where the last equality follows combining  \Cref{prop:etalepullbackofV} and \Cref{Prop:MHS on Cl are constant}.  
\end{proof}
Let $\lambda$ be a partition of $n$. Define the graded mixed Hodge structure 
\[ 
H_{\lambda} \defeq  \bigoplus_{j \in \mathbb{Z}} H_{\lambda,j}[-j].
\] 
Then \Cref{Prop:ICs as equi sheaves} implies that
\begin{align}\label{eqn: H in the decompo}
     \mathcal{F}_{n,\SE} \simeq \bigoplus_{\substack{\lambda \vdash n}} \,\left(q_{\lambda\ast} H_{\lambda}\right)^{G_{\lambda}}.
\end{align}
\begin{proposition} \label{prop:factorizingH}
Let  $\lambda=(\lambda_1,\ldots, \lambda_{l(\lambda)})$ be a partition of $n$. Then there is  an isomorphism of $G_{\lambda}$-equivariant mixed Hodge structures
\[ 
H_{\lambda} \simeq \bigotimes_{i=1}^{l(\lambda)} H_{(\lambda_i)},
\] 
where the group $G_{\lambda}$ acts on $ \bigotimes_{i=1}^{l(\lambda)} H_{(\lambda_i)}$ by permuting the  factors and $(\lambda_i)$ denotes the partition of $\lambda_i$ with a single part of size $\lambda_i$. 
\end{proposition}

\begin{proof}
Set $l \defeq l(\lambda)$.
By \Cref{lemma:addition-is-etale}, the addition of cycles map $\cup_\lambda \colon A_{\lambda}\to \Sym^n(S)$ is \'etale.  By \Cref{prop:sheaffactorization} we have that 
\begin{equation}
\label{eq:F,n,V as box product}
\cup_\lambda^* \mathcal{F}_{n,\SE} \simeq \left.\left(\bigboxtimes_{i=1}^{l} \mathcal{F}_{\lambda_i,\SE}\right)\right|_{ A_{\lambda}}.
\end{equation}
Denote by 
\[
\begin{tikzcd}
\Delta^\lambda_S\colon C_\lambda(S) \arrow[hook]{r} & A_{\lambda}
\end{tikzcd}
\]
the diagonal embedding $(p_1,\ldots, p_l) \mapsto (\lambda_1 [p_1],\ldots, \lambda_l [p_l])$.
   We now analyze the supports of the summands on both sides of Equation \eqref{eq:F,n,V as box product}.  On  the left, we have that
    \[
    \cup_\lambda^* \mathcal{F}_{n,\SE}=\bigoplus_{\substack{\mu \vdash n \\ j \in \mathbb{Z}}} \cup_\lambda^*\mathrm{IC}_{\overline{\Sym^n_{\lambda}(S)}}(\mathbb{V}_{\mu,j})[-j]   = \left(\left(\Delta^{\lambda}_S\right)_{\ast}H_{\lambda}\right) \oplus E,
    \]
    where $E$ is a direct sum of shifts of simple mixed Hodge modules, with supports different from $C_\lambda(S)$.
On the other hand, set $R_\lambda= \bigotimes_{i=1}^{l(\lambda)} H_{(\lambda_i)}$.
Then  the right hand side of Equation \eqref{eq:F,n,V as box product} decomposes as
    \[ 
\left(\left(\Delta^{\lambda}_S\right)_{\ast}R_{\lambda}\right) \oplus F,
    \] 
    where $F$ can be written as a direct sum of simple mixed Hodge modules with supports different from the copy of $C_\lambda(S)$ in $A_{\lambda}$. Thus, by comparing supports we obtain the isomorphism
    \[ H_{\lambda} \simeq R_\lambda = \bigotimes_{i=1}^{l(\lambda)} H_{(\lambda_i)},
    \]
    as required.
\end{proof}

As an application of  \Cref{prop:factorizingH} and \Cref{Prop:ICs as equi sheaves}, we express the mixed Hodge modules $\CF_{n, \SE}$ in terms of the mixed Hodge structures $H_{(m)}$ via a plethystic expression and compute the cohomology of $\Quot^n(\SE)$. 

\begin{theorem}
\label{thm:sheaf and coho as Sym}
Let $S$ be a smooth quasiprojective surface and $\SE$ a locally free sheaf on $S$. Then there is an isomorphism
\begin{align}\label{eqn: fighissimo}
    \bigoplus_{n \geq 0} \mathcal{F}_{n,\SE} \simeq \Sym_{\cup}\left( \bigoplus_{n \geq 1} \left(q_{(n)\ast}\mathbb{Q}_{S}\otimes H_{(n)}\right)\right)\,\in\,\widehat{\derived}(\Sym(S)).
\end{align}
In particular, $\CF_{n, \SE}$ depends only on $\rk \SE$.
Moreover, we have an isomorphism of $\BZ$-graded mixed Hodge structures
\begin{equation}
\label{eq:Coho of Quot Sym expression}
\bigoplus_{n \geq 0} \HH^{\ast}_{(c)}(\Quot^n(\SE)) \simeq \Sym \left( \bigoplus_{n \geq 1} \HH^{\ast}_{(c)}(S) \otimes H_{(n)}\right).
\end{equation}
\end{theorem}
\begin{proof}
Let $l\geq 0$ and  denote by $\cup^l\colon \Sym(S)^{l} \to \Sym(S)$ the $l$-fold addition map. 
By the projection formula, we have that 
\begin{align}\label{eqn: prok in proof}
 &q_{(n)\ast} H_{(n)}= q_{(n)\ast}\mathbb{Q}_{S}\otimes H_{(n)}.
\end{align}
By definition of $\Sym_{\cup}$ (cf.~\Cref{subsec: Plethystic operations}), we have that
\begin{align*}
    \Sym_{\cup}\left( \bigoplus_{n \geq 1} q_{(n)\ast}H_{(n)}\right)&=\bigoplus_{l\geq 0}\Sym_{\cup}^{l}\left( \bigoplus_{n \geq 1} q_{(n)\ast}H_{(n)}\right)\\
    &=\bigoplus_{l\geq 0}\left(\cup^l_!\left(\bigoplus_{n \geq 1} q_{(n)\ast}H_{(n)}\right)^{\boxtimes l}\right)^{\FS_l}\\
&=\bigoplus_{l\geq 0}\left(\cup^l_!\left(\bigoplus_{(k_{1},\ldots,k_{l})}\bigboxtimes_{1\leq i\leq l}q_{(k_i)\ast}H_{(k_i)}
\right)\right)^{\FS_l}\\
&=\bigoplus_{l\geq 0}\left(\bigoplus_{\lambda_{1}\geq \cdots \geq\lambda_{l}}\left(\cup^l_!\left(\bigboxtimes_{1\leq i\leq l}q_{(\lambda_i)\ast}H_{(\lambda_i)}
\right)\right)^{G_{\lambda}}\right)\\
&= \bigoplus_{n\geq 0}\left(\bigoplus_{\lambda\vdash n}\left({q_{\lambda}}_{*}H_{\lambda}\right)^{G_{\lambda}}\right),
\end{align*}
where in the last line we applied \Cref{prop:factorizingH} and the relative Künneth formula. The last expression is equal to $\bigoplus_{n \geq 0} \mathcal{F}_{n,\SE}$ by \Cref{Prop:ICs as equi sheaves}, see also \eqref{eqn: H in the decompo}. Combining this identification with \eqref{eqn: prok in proof} proves the identity \eqref{eqn: fighissimo}.

Finally, the isomorphisms of mixed Hodge structures \eqref{eq:Coho of Quot Sym expression} are  obtained from \eqref{eqn: fighissimo} by taking (derived) global sections and (derived) compactly supported  global sections, respectively. 
\end{proof}

We finally compute  the mixed Hodge structure $H_{(n)}$ appearing in \Cref{thm:sheaf and coho as Sym}, concluding the proof of \Cref{thm: Computation of decomposition}.

\begin{proposition}
\label{prop:computation of Hn}
Let $n\geq 1$ and set $r=\rk \SE$. Then there is an 
isomorphism of mixed Hodge structures
    \[ 
    H_{(n)} \simeq \HH^{\ast}(\mathbb{P}^{r-1}) \otimes \mathbb{L}^{(n-1)r}.
    \]
\end{proposition}
\begin{proof}
By \Cref{Prop:ICs as equi sheaves} and \Cref{prop:factorizingH} we have that
\begin{equation}
\label{eqn:Fnvdecomposition}
    \mathcal{F}_{n,\SE} = \bigoplus_{\lambda \vdash n} \left(q_{\lambda\ast} H_{\lambda}\right)^{G_{\lambda}}=\bigoplus_{\lambda \vdash n} q_{\lambda *}\left(\bigotimes_{i=1}^{l(\lambda)} H_{(\lambda_i)}\right)^{G_{\lambda}}.
\end{equation}
Let $p\in S$ be any point and consider the closed immersion $\iota\colon [n] p\hookrightarrow \Sym^{n}(S)$. By proper base change and \eqref{eqn: def Fn}, we have that
\begin{align}\label{eqn: proof eq 1}
    \iota^*  \mathcal{F}_{n,\SE} \simeq \HH^{\ast}(\Quot_{n,r}^{0}), 
\end{align}
where we used that $\pi_{n,\SE}^{-1}(n[p])$ is isomorphic to 
$\Quot_{n,r}^0$. On the other hand, applying proper base change to \eqref{eqn:Fnvdecomposition} yields
\begin{align}\label{eqn_ proof eqn 2}
     \iota^*  \mathcal{F}_{n,\SE} \simeq \bigoplus_{\lambda \vdash n} \left(\bigotimes_{i=1}^{l(\lambda)} H_{(\lambda_i)}\right)^{G_{\lambda}}.
\end{align}
Combining \eqref{eqn: proof eq 1} with \eqref{eqn_ proof eqn 2} yields the identification
\begin{equation}\label{eqn:cohomologyofquot}
\HH^{\ast}(\Quot_{n,r}^{0}) \simeq \bigoplus_{\lambda \vdash n}\left(\bigotimes_{i=1}^{l(\lambda)} H_{(\lambda_i)}\right)^{G_{\lambda}}.
\end{equation} 
By \cite[Cor.~3.10]{NY_lectures}, the punctual Quot scheme $\Quot_{n,r}^{0}$ admits a cellular decomposition; in particular, its cohomology $\HH^{\ast}(\Quot_{n,r}^{0}) $ is pure and satisfies
\begin{align}\label{eqn: sym for QUot zzzz}
    \bigoplus_{n \geq 0} \HH^{\ast}(\Quot_{n,r}^{0}) \simeq \Sym \left( \bigoplus_{n \geq 1} \HH^{\ast}(\mathbb{P}^{r-1}) \otimes \mathbb{L}^{(n-1)r}\right).
\end{align}
On the other hand, \eqref{eqn:cohomologyofquot} and the same manipulations as in the proof of \Cref{thm:sheaf and coho as Sym} imply that 
\begin{align}\label{eqn: last Sym QUot 0}
   \bigoplus_{n \geq 0} \HH^{\ast}(\Quot_{n,r}^{0}) \simeq \Sym \left( \bigoplus_{n \geq 1} H_{(n)}\right). 
\end{align}
The result then follows combining \eqref{eqn: sym for QUot zzzz} and \eqref{eqn: last Sym QUot 0}.
\end{proof}

The following is a straightforward consequence of \Cref{prop:computation of Hn} and \Cref{thm:sheaf and coho as Sym}, and gives a positive answer to Question \ref{Question: Mozgovoy}, posed by Mozgovoy \cite{mozgovoy2019motivic}.

\begin{corollary}
\label{cor: Purity of Quot}
 Let $S$ be a smooth quasiprojective surface such that $\HH^{\ast}(S)$, resp.~$\HH^{\ast}_{c}(S)$, is a pure mixed Hodge structure of weight $0$. Then, for all $n\geq 0$, we have that $\HH^{\ast}(\Quot^n(\SE))$, resp.~$\HH_{c}^{*}(\Quot^n(\SE))$, is a pure mixed Hodge structure of weight $0$, i.e.~$\HH^k_{(c)}(\Quot^n(\SE))$ is a pure mixed Hodge structure of weight $k$ for all $k \geq 0$.  
\end{corollary}

\begin{example}
\label{example: Thm B generalizes Gottsche Soeergel}
Combining \Cref{prop:factorizingH} and \Cref{prop:factorizingH} yields, for all $r\geq 1$,
\begin{equation}
H_{\lambda} \simeq \bigotimes_{i=1}^{l(\lambda)}H_{(\lambda_i)} = \bigotimes_{i=1}^{l(\lambda)} \left(\HH^{\ast}(\mathbb{P}^{r-1}) \otimes \mathbb{L}^{(\lambda_i-1)r} \right) \simeq \HH^{\ast}(\mathbb{P}^{r-1})^{l(\lambda)} \otimes \mathbb{L}^{(n-l(\lambda))r}.
\end{equation} 
When $r=1$, i.e.~$\Quot^n(\SE) \simeq \mathrm{Hilb}^n(S)$, we have 
\[ 
H_{\lambda} \simeq \mathbb{L}^{n-l(\lambda)},
\] 
endowed with the  trivial action of $G_\lambda$. Thus $H_{\lambda,j} = \HH^{j-2l(\lambda)}( H_{\lambda}) = 0$ for all $j \neq 2n$. Therefore by 
\Cref{prop:etalepullbackofV} we have
\[\mathbb{V}_{\lambda,2n}[2l(\lambda)-2n] \simeq \mathbb{L}^{n-l(\lambda)}\]
and $\mathbb{V}_{\lambda,j}=0$ for  $j\neq 2n$. Thus we get 
\[ 
\mathcal{F}_{n,\OO_S} \simeq \bigoplus_{\lambda \vdash n} \mathrm{IC}_{\overline{\Sym^n_{\lambda}(S)}}(\mathbb{V}_{\lambda, 2n})[-2n] = \bigoplus_{\lambda \vdash n}\mathrm{IC}_{\overline{\Sym^n_{\lambda}(S)}}\left(\mathbb{L}^{n-l(\lambda)}\right)[-2l(\lambda)],
\] which recovers \cite[Thm.~4]{gotsoergelhilbert}.
\end{example}

\section{Applications}\label{section:applications}
In this section, we will derive applications of our main results.

\subsection{Hodge--Deligne series and factorization into rank 1 series}
\label{sec:HD-series}
In this section, we compute the Hodge--Deligne series for $\Quot^n(\SE)$, both ordinary and compactly supported. Before recalling its definition, we set up the plethystic formalism useful to derive the compact formula anticipated in \Cref{thm: intro HD}. Set $R = \BZ[u^{\pm 1},v^{\pm 1},z^{\pm 1}]$. Consider the (invertible) group homomorphism
\[
\map{\Sym}{(tR \llbracket t \rrbracket,+)}{(1+tR \llbracket t \rrbracket,\cdot),}
\]
defined as the plethystic exponential attached to the Adams operators $\psi_e^{\mathrm{super}}(\boldit m) = (-1)^{(e-1)k}\boldit m^{e}$, where $\boldit m = u^pv^qz^kt^m$ is a generic monomial appearing in a power series 
\[
f = \sum_{m\geq 1}\sum_{p,q,k}a_{m,p,q,k}u^pv^qz^kt^m \in tR \llbracket t \rrbracket.
\]
The sign $(-1)^{(e-1)k}$ reflects the super sign rule $x \otimes y \mapsto (-1)^{\lvert x\rvert \lvert y \rvert}y \otimes x$ in the category of cohomologically graded vector spaces, where the cohomological degree is recorded by the variable $z$. Explicitly, $\Sym$ acts on a basic monomial as above as
\[
\Sym(a\boldit m) = 
\begin{cases}
    (1-\boldit m)^{-a} & \text{if }k\text{ is even,} \\
    (1+\boldit m)^{a} & \text{if }k\text{ is odd.}
\end{cases}
\]
Consequently, for a power series as above, one has
\begin{equation}
\label{general-Sym-formula}
\Sym \left(\sum_{m\geq 1}\sum_{p,q,k}a_{m,p,q,k}u^pv^qz^kt^m \right) = \prod_{m,p,q,k} (1-(-1)^ku^pv^qz^kt^m)^{(-1)^{k+1}a_{m,p,q,k}}.
\end{equation}
This expression can also be written as
\[
\exp\left(\sum_{e\geq 1}\frac{1}{e}\psi_e^{\mathrm{super}}(f) \right).
\]
With the plethystic formalism now in place, recall that the Hodge--Deligne series of a cohomologically graded mixed Hodge structure $\mathscr L$ is the formal power series
\begin{align}\label{eqn: HD graaaa}
    \mathsf h(\mathscr L;x,y,z) = \sum_{p,q,k\in \BZ}\dim(\Gr^p_H(\Gr^W_{p+q}(\mathrm H^k(\mathscr L))))u^pv^qz^k.
\end{align}
If $Y$ is a finite type separated $\BC$-scheme, we define $\mathsf h_{(c)}(Y;u,v,z)$ to be the Hodge--Deligne series attached to $\HH_{(c)}^\ast(Y,\BQ)$. We now compute the Hodge--Deligne series, both ordinary and compactly supported, for all Quot schemes at once.
\smallbreak
Let us set $\mathsf q \defeq \mathsf h_c(\BA^1;u,v,z) = uvz^2$, so that $\mathsf h(\BP^{r-1};u,v,z) = \sum_{0\leq j\leq r-1}\mathsf q^j$, and denote by 
\begin{align*}
    \mathsf{HD}^{(c)}_r(t)=\sum_{n\geq 0}\mathsf h_{(c)}(\Quot^n(\SE);u,v,z)t^n
\end{align*}
the  generating series of Hodge--Deligne series,  where $S$ is a smooth quasiprojective surface and $\SE$ is a locally free sheaf on $S$.
\begin{theorem}
\label{thm: HD poly}
Let $S$ be a smooth quasiprojective surface with Hodge--Deligne polynomial 
\[
\mathsf h_{(c)}(S;u,v,z) = \sum_{p,q,k}h_{p,q,{(c)}}^ku^pv^qz^k,
\]
where $h_{p,q,(c)}^k = h^{p,q}(\HH^k_{(c)}(S,\BQ))$. Let $\SE$ be a rank $r$ locally free sheaf on $S$. There are identities
\begin{align*}
  \mathsf{HD}^{(c)}_r(t)&= \Sym \left(\mathsf h_{(c)}(S;u,v,z)\frac{1-\mathsf q^r}{1-\mathsf q}\frac{t}{1-t\mathsf q^r}\right) \\
&=\prod_{m\geq 1}\prod_{j=0}^{ r-1}\prod_{p,q,k}\left(1-(-1)^kt^mu^{p+j+r(m-1)} v^{q+j+r(m-1)}z^{k+2j+2r(m-1)}\right)^{(-1)^{k+1}h_{p,q,(c)}^k}.    
\end{align*}
\end{theorem}

\begin{proof}
By multiplicativity of $\mathsf h_{(c)}$, we have
\[
\mathsf h_{(c)}(\HH^\ast(S \times \BP^{r-1}) \otimes \BL^{r(m-1)};u,v,z) = \mathsf h_{(c)}(S;u,v,z)\mathsf q^{r(m-1)}\sum_{0\leq j\leq r-1}\mathsf q^j.
\]
Summing over $m\geq 1$ yields
\[
\sum_{m\geq 1} \mathsf h_{(c)}(\HH^\ast(S \times \BP^{r-1})\otimes \BL^{r(m-1)};u,v,z) t^m = \mathsf h_{(c)}(S;u,v,z)\frac{1-\mathsf q^r}{1-\mathsf q}\frac{t}{1-t\mathsf q^r}.
\]
\Cref{thm:sheaf and coho as Sym} then yield identities 
\begin{equation}
\label{Sym-identity-Hodge-deligne}
\sum_{n\geq 0}\,\mathsf h_{(c)}(\Quot^n(\SE);u,v,z)t^n = \Sym \left(\mathsf h_{(c)}(S;u,v,z)\frac{1-\mathsf q^r}{1-\mathsf q}\frac{t}{1-t\mathsf q^r}\right).
\end{equation}
So the claimed infinite product identities follow formally from \Cref{general-Sym-formula}.
\end{proof}
For a separated, finite type $\BC$-scheme $Y$, define the Poincar\'e polynomial $P(Y,z)$ and the compactly supported Poincar\'e polynomial $P_c(Y,z)$ by
\[
P_{(c)}(Y,z) = \sum_{i\geq 0} b_i^{(c)}z^i,
\]
where $b_i = \dim_{\BQ} \HH^i(Y(\BC),\BQ)$ and $b_i^c = \dim_{\BQ} \HH^i_c(Y(\BC),\BQ)$ are the Betti numbers.  Of course, if $S$ is a projective surface, the two expressions agree. 
Specialising $u=v=1$ in \Cref{thm: HD poly}, we obtain a closed formula for the Poincaré polynomial of $\Quot^n(\SE)$, generalising the classical formula of Göttsche \cite{Gott2}.
\begin{corollary}\label{cor: Poincare}
Let $S$ be a smooth quasiprojective surface with Betti numbers $b_i^{(c)}$, for $i=0,\ldots,4$. Let $\SE$ be a locally free sheaf of rank $r$ over $S$. There are identities
\begin{align}
\label{eqn: Poincare for lemma}
    \sum_{n\geq 0}P_{(c)}(\Quot^n(\SE),z)t^n=
\prod_{m\geq 1}\prod_{j=0}^{r-1}\prod_{i=0}^4 \left(1-(-1)^it^m z^{2j+2r(m-1)+i}\right)^{(-1)^{i+1}b_i^{(c)}}.
\end{align}
\end{corollary}
We remark that \Cref{cor: Poincare} computes the Poincar\'e polynomials of $ \Quot^n(\SE)$ also when $S$ is a quasiprojective surface with possibly non-pure cohomology. In particular, this implies that \eqref{eqn: Poincare for lemma} is not a motivic specialization of the motivic series computed in \cite{mozgovoy2019motivic}. 

\subsubsection{Rank $r$-to-1 factorization}
Set $R=\BZ[u^{\pm 1},v^{\pm 1},z^{\pm 1}]$. Denote by 
\[
\map{\Log_{\mathrm{super}}}{1+tR\llbracket t \rrbracket}{tR\llbracket t \rrbracket}
\]
the inverse of the $\Sym$ operator.
For an element $a \in R$ and a power series $g \in 1+tR\llbracket t \rrbracket$, define
\[
g^{a} \defeq \Sym\left( a\Log_{\mathrm{super}}(g)\right) \in 1+tR\llbracket t \rrbracket.
\]
This operation satisfies 
\[
\Sym(ag) = \Sym(g)^a,
\]
and in particular defines a power structure on $R$. With this formalism in place, we are able to factor the Hodge--Deligne series in rank $r$, cf.~\Cref{thm: HD poly}, as a product of $r$ rank $1$ Hodge--Deligne series, as the next result illustrates.

\begin{proposition}\label{prop: rank r to 1}
Let $S$ be a smooth quasiprojective surface and $\SE$ a rank $r$ locally free sheaf on $S$.  Then, there is an identity
\[
\mathsf{HD}^{(c)}_r(t) = \prod_{j=0}^{r-1} \mathsf{HD}^{(c)}_1(t\mathsf q^{r-1})^{\mathsf q^{j+1-r}}\,\in\,1+tR\llbracket t \rrbracket.
\]
\end{proposition}

\begin{proof}
Our starting point is \Cref{Sym-identity-Hodge-deligne}. A straightforward verification shows that, if we set $\widetilde t = t\mathsf q^{r-1}$, we have
\[
\frac{1-\mathsf q^r}{1-\mathsf q}\frac{t}{1-t\mathsf q^r} = \left(\sum_{j=0}^{r-1}\,\mathsf q^{j+1-r}\right)\frac{\widetilde t}{1-\widetilde t\mathsf q}.
\]
Therefore, \Cref{Sym-identity-Hodge-deligne} becomes
\begin{align*}
\mathsf{HD}^{(c)}_r(t) 
&= \Sym \left(\sum_{j=0}^{r-1}\,\mathsf q^{j+1-r}\mathsf h_{(c)}(S;u,v,z) \frac{t\mathsf q^{r-1}}{1-t\mathsf q^{r-1}\cdot \mathsf q}\right) \\
&= \prod_{j=0}^{r-1}\Sym\left(\mathsf q^{j+1-r}\mathsf h_{(c)}(S;u,v,z)\frac{t\mathsf q^{r-1}}{1-t\mathsf q^{r-1}\cdot \mathsf q}\right) \\
&= \prod_{j=0}^{r-1}\Sym\left(\mathsf h_{(c)}(S;u,v,z)\frac{t\mathsf q^{r-1}}{1-t\mathsf q^{r-1}\cdot \mathsf q}\right)^{\mathsf q^{j+1-r}} \\
&=\prod_{j=0}^{r-1} \mathsf{HD}^{(c)}_1(t\mathsf q^{r-1})^{\mathsf q^{j+1-r}},
\end{align*}
as required.
\end{proof}

\subsection{Kirwan surjectivity for Quot schemes}\label{subsec:Kirwan Surjectivity}
We follow the notation of  \Cref{subsection: KIrwan Surjectivity intro}. Consider the commutative diagram
\[
\begin{tikzcd}
    \Quot^n(\SE)\arrow[r,hook, "j"]\arrow[swap,bend right=20]{rdd}{\pi_{n,\SE}} &  \mathrm{Fr}_{n}(\SE)\arrow[d, "p"]\arrow[dd, bend left=70, "q"]\\
        & \CCoh^{n}(S)\arrow[d, "\sigma_n"]\\
        & \Sym^{n}(S),
    \end{tikzcd}
\]
where $j$ is the open immersion of the semistable locus. We let $i\colon Z_{n,\SE} \hookrightarrow \mathrm{Fr}_{n}(\SE)$ be the closed complement of 
$\Quot^n(\SE)$ in $\mathrm{Fr}_{n}(\SE)$. 
Kirwan surjectivity for $j$ is a corollary of \Cref{thm:Main thm inro}.
\begin{corollary}
\label{cor:Kirwan Surjectivity}
The unit $\id \to j_*j^*$ induces a split surjection of pure complexes of mixed Hodge modules 
\begin{equation}\label{eqn: duals in proof}
\begin{tikzcd}
q_* \mathbb{D}\mathbb{Q}_{\mathrm{Fr}_n(\SE)} \arrow{r} & \pi_{n,\SE\ast}\mathbb{D}\mathbb{Q}_{\Quot^n(\SE)},
\end{tikzcd}
\end{equation}
where $\mathbb{D}\mathbb{Q}$ denotes the dualising mixed Hodge module.

Similarly, the open pullback in Borel--Moore homology induces a surjection of mixed Hodge structures
\[
\begin{tikzcd}
\HH_{\ast}^{\mathrm{BM}}(\CCoh^{n}(S))\otimes\mathbb{L}^{-rn}\cong \HH_{\ast}^{\mathrm{BM}}(\mathrm{Fr}_{n}(\SE))\arrow[two heads]{r} & \HH_{\ast}^{\mathrm{BM}}(\Quot^n(\SE)).
\end{tikzcd}
\]
Dually, there is an injection in compactly supported cohomology
    \[
\begin{tikzcd}
    \HH^{\ast}_{c}(\Quot^n(\SE))\arrow[hook]{r} &  \HH^{\ast}_{c}(\mathrm{Fr}_{n}(\SE))\cong \HH^{\ast}_{c}(\CCoh^{n}(S))\otimes \mathbb{L}^{rn}.
\end{tikzcd}
    \]
\end{corollary}
\begin{proof}
Consider the localization exact triangle
\[
\begin{tikzcd}
i_{*}\mathbb{D}\mathbb{Q}_{Z_{n,\SE}}\arrow{r} & \mathbb{D}\mathbb{Q}_{\mathrm{Fr}_{n}(\SE)}\arrow{r} & j_{*}\mathbb{D}\mathbb{Q}_{\Quot^n(\SE)}\arrow{r} & i_{*}\mathbb{D}\mathbb{Q}_{Z_{n,\SE}}[1]
\end{tikzcd}
\]
on $\mathrm{Fr}_{n}(\SE)$.
Applying $q_{*}$ yields the exact triangle
\begin{equation}
\label{eq: triangle in kirwan}
\begin{tikzcd}
(q\circ i)_{*}\mathbb{D}\mathbb{Q}_{Z_{n,\SE}}\arrow{r} & q_{*}\mathbb{D}\mathbb{Q}_{\mathrm{Fr}_{n}(\SE)}\arrow{r} & {\pi_{n,\SE\ast}}\mathbb{D}\mathbb{Q}_{\Quot^n(\SE)}\arrow{r}{\delta} & (q\circ i)_{*}\mathbb{D}\mathbb{Q}_{Z_{n,\SE}}[1].
\end{tikzcd}
\end{equation}
Since exceptional pullbacks and usual pushforwards do not decrease weights, we obtain that the weights of $(q\circ i)_{*}\mathbb{D}\mathbb{Q}_{Z_{n,\SE}}$ are greater than or equal to $0$ and hence the weights of $(q\circ i)_{*}\mathbb{D}\mathbb{Q}_{Z_{n,\SE}}[1]$ are $\geq 1$. Since the Quot-to-Chow map $\pi_{n,\SE}$ is proper, \Cref{thm:Main thm inro}  implies that ${\pi_{n,\SE\ast}}\mathbb{D}\mathbb{Q}_{\Quot^n(\SE)}$ is a pure mixed Hodge module of weight $0$. Therefore the map $\delta$ vanishes and hence the triangle \eqref{eq: triangle in kirwan} splits. 
Pushing \eqref{eqn: duals in proof} to a point yields the required  surjection in Borel--Moore homology, and dualizing yields the required injection in compactly supported cohomology.
\end{proof}

\begin{remark}
\label{rem:Kirwan surjectivity lit}
The stack $\mathrm{Fr}_n(\SE)$ has a natural enhancement to a \emph{quasi-smooth} derived stack, and the inclusion $\Quot^n(\SE) \subset \mathrm{Fr}_{n}(\SE)$ can be seen as the semistable locus of a natural (derived) $\Theta$-stratification on $\mathrm{Fr}_n(\SE)$. In this context, 
    Kirwan surjectivity is known for some large classes of $\Theta$-stratified quasi-smooth derived stacks. Some relevant examples are given by Hennecart--Kinjo \cite[Thm.~6.6]{arXiv:2509.21298} in the case of symplectic stacks 
    and by Halpern-Leistner \cite[Thm.~3.1]{zbMATH07272594} in the case of quasi-smooth stacks satisfying a technical condition $(\dagger)$, which is satisfied by symplectic stacks. However, the derived enhancement of $\mathrm{Fr}_n(\SE)$  does not satisfy the condition $(\dagger)$ (in particular, it is not symplectic) hence one cannot directly deduce Kirwan surjectivity for  $\Quot^n(\SE) \subset \mathrm{Fr}_{n}(\SE)$ from their results.
\end{remark}
    
\subsection{Purity for the stack of 0-dimensional coherent sheaves}\label{sec: purity Coh}
In this section we give a different proof of the following purity result of Davison, which is \cite[Thm.~A]{Davisonpreprojective} applied to the Jordan quiver, see also \cite[Prop.~9.1]{BDHSM}.

\begin{corollary}\label{cor:Purity of Coh}
Let $S$ be a smooth quasiprojective surface such that $\HH_{\ast}^{\mathrm{BM}}(S)$ is a pure mixed Hodge structure of weight $0$. Then, for all $n\geq 0$, we have that $\HH_{\ast}^{\mathrm{BM}}(\CCoh^{n}(S))$ is a pure mixed Hodge structure of weight $0$.
\end{corollary}

\Cref{cor:Purity of Coh} will  follow immediately from  \Cref{cor: Purity of Quot} and the following stabilization result. In the following, by convention we consider the Borel--Moore homology $\HH_{i}^{\mathrm{BM}}(-)$ as a complex supported in homological degree $i$. 

\begin{proposition}
    \label{prop: Homology of Quot stabilizes to that of Coh}
Let $S$ be a smooth quasiprojective surface. Given $k\in\mathbb{Z}$ and $n\geq 0$, there exists $r_{0}\geq n$ such that for all  $r>r_{0}$ we have an isomorphism
\[
\begin{tikzcd}
j^*\circ p^!\colon \HH_{k}^{\mathrm{BM}}(\CCoh^{n}(S))\otimes\mathbb{L}^{-rn}\arrow{r}{\sim} & \HH_{k+2rn}^{\mathrm{BM}}(\Quot^n(\OO^{\oplus r}_S)).
\end{tikzcd}
\]
In particular, we have
\[
\HH_{*}^{\mathrm{BM}}(\CCoh^{n}(S))\cong \varprojlim_r\HH_{*}^{\mathrm{BM}}(\Quot^n(\OO^{\oplus r}_S))\otimes \mathbb{L}^{rn}.
\]
\end{proposition}

\begin{proof}
   Assume that $r\geq n$ and define $U_{n,r}$ to be the open subscheme of $\Quot^n(\OO^{\oplus r}_S)$ consisting of  quotients $\phi\colon \OO_S^{\oplus r}\onto \SQ$ such that the induced map 
    \[
    \begin{tikzcd}
     \BC^r\arrow[hook]{r} & \HH^{0}(S,\OO_S^{\oplus r})\arrow{r} & \HH^{0}(S,\SQ)\cong \mathbb{C}^{n} 
    \end{tikzcd}
    \]
    is a surjection. 
    Here, the map $\mathbb{C}^{r}\into \HH^{0}(S,\OO_S^{\oplus r})$ is the standard $\mathbb{C}$-linear inclusion of $\mathbb{C}^{r}$, coordinate-wise.
    Let $Y_{n,r}$ be the closed complement of $U_{n,r}$ in $\mathrm{Fr}_{n}(\OO_S^{\oplus r})$. Clearly, we have an inclusion $Z_{n,\OO_S^{\oplus r}}\subset Y_{n,r}$. We claim that $Y_{n,r}$ has codimension $r-n+1$ in $\mathrm{Fr}_{n}(\OO_S^{\oplus r})$. Indeed, every fibre of the map $p\colon \mathrm{Fr}_{n}(\OO_S^{\oplus r}) \to \CCoh^n(S)$ can be identified with the space of linear maps $\mathrm{Hom}(\mathbb{C}^{r},\mathbb{C}^{n})$ and every fibre of the induced map 
\[
\begin{tikzcd}
Y_{n,r}\arrow[hook]{r} & \mathrm{Fr}_{n}(\OO_S^{\oplus r})\arrow{r}{p} & \CCoh^{n}(S)
\end{tikzcd}
\]
can then be identified with the closed subspace of non-surjective linear maps in $\mathrm{Hom}(\mathbb{C}^{r},\mathbb{C}^{n})$. One verifies that this space has codimension $r-n+1$ in $\mathrm{Hom}(\mathbb{C}^{r},\mathbb{C}^{n})$ and hence $Y_{n,r}$ has codimension $r-n+1$ in $\mathrm{Fr}_{n}(\OO^{\oplus r})$. Since $\mathrm{dim}\,\mathrm{Fr}_{n}(\OO_S^{\oplus r})= \mathrm{dim}\,\Quot^n(\OO^{\oplus r}_S)=rn+n$, we obtain that 
\begin{equation}\label{eq: upper bound of dim of Z}
\mathrm{dim}\,Z_{n,\OO^{\oplus r}}\leq \mathrm{dim}\,Y_{n,r} = rn - r+2n-1.
\end{equation}
We now consider the localization long exact sequence in Borel--Moore homology
\[
\cdots \rightarrow \HH_{k+2rn}^{\mathrm{BM}}(Z_{n,\OO^{\oplus r}})\rightarrow \HH_{k+2rn}^{\mathrm{BM}}(\mathrm{Fr}_{n}(\OO^{\oplus r}))\rightarrow \HH_{k+2rn}^{\mathrm{BM}}(\Quot^n(\OO^{\oplus r}_S))\rightarrow \HH_{k+2rn-1}^{\mathrm{BM}}(Z_{n,\OO^{\oplus r}})\rightarrow \cdots
\]
and we note that if $r$ is such that $2r>4n-k$, then Equation \eqref{eq: upper bound of dim of Z} combined with the long exact sequence yields an isomorphism
\[
\HH_{k+2rn}^{\mathrm{BM}}(\mathrm{Fr}_{n}(\OO_S^{\oplus r}))\cong \HH_{k+2rn}^{\mathrm{BM}}(\Quot^n(\OO^{\oplus r}_S)).
\]
Along with the isomorphism $\HH_{\ast}^{\mathrm{BM}}(\mathrm{Fr}_{n}(\OO_S^{\oplus r}))\cong \HH_{\ast}^{\mathrm{BM}}(\CCoh^{n}(S))\otimes\mathbb{L}^{-rn}$, this finishes the proof.
\end{proof}

\begin{remark}\label{ThmAisEquivToKirwanandDavison}
In view of the results of Sections \ref{subsec:Kirwan Surjectivity} and \ref{sec: purity Coh}, consider the following statements:
\begin{enumerate}
\item [\mylabel{cond-A}{(A)}] \cref{thm:Main thm inro} holds for all locally free sheaves $\SE$ on $S$,
\item [\mylabel{cond-B}{(B)}] The natural map $q_{*}\mathbb{D}\mathbb{Q}_{\mathrm{Fr}_{n}(\SE)}\rightarrow {\pi_{n,\SE\ast}}\mathbb{D}\mathbb{Q}_{\Quot^n(\SE)}$ is a split epimorphism for all locally free sheaves  $\SE$ on $S$,
\item [\mylabel{cond-C}{(C)}] The pushforward $\sigma_{n\ast}\mathbb{DQ}_{\CCoh^n(S)}$ is  pure of weight 0.
\end{enumerate}
Then, combining Corollaries \ref{cor:Kirwan Surjectivity} and \ref{cor:Purity of Coh} with \Cref{thm:Main thm inro} yields the following equivalence:
\begin{align}\label{eqn: equivalence}
\textrm{\ref{cond-A} holds $\iff$ \ref{cond-B} and \ref{cond-C} hold.}    
\end{align}
The proof of \eqref{eqn: equivalence} follows by similar arguments to the ones used to prove Corollaries \ref{cor:Kirwan Surjectivity} and \ref{cor:Purity of Coh}. Notice also that \ref{cond-C} does not  imply \cref{cor:Purity of Coh} when $S$ is not projective.

\end{remark}

\bibliographystyle{amsplain-nodash}
\bibliography{bib}

\bigskip
\noindent
\emph{Shivang Jindal}, \texttt{shivang.jindal@epfl.ch} \\
\textsc{EPFL, Institute of Mathematics, CH-1015 Lausanne,
Switzerland}

\medskip
\noindent
\emph{Archi Kaushik}, \texttt{archi.kaushik@epfl.ch} \\
\textsc{EPFL, Institute of Mathematics, CH-1015 Lausanne,
Switzerland}

\medskip
\noindent
\emph{Sergej Monavari}, \texttt{sergej.monavari@unipd.it} \\
\textsc{Università di Padova, Via Trieste 63, 35121 Padova, Italy}

\medskip
\noindent
\emph{Andrea T. Ricolfi}, 
\texttt{aricolfi@sissa.it} \\
\textsc{SISSA, Via Bonomea 265, 34136 Trieste, Italy}

\medskip
\noindent
\emph{Sebastian Schlegel Mejia}, \texttt{sebastian.schlegelmejia@epfl.ch} \\
\textsc{EPFL, Institute of Mathematics, CH-1015 Lausanne,
Switzerland}

\end{document}